\documentclass[journal,10pt]{arXiv}

\usepackage[utf8]{inputenc}

\usepackage{amsmath}
\usepackage{amsthm}
\usepackage{amssymb}
\usepackage{graphicx}
\graphicspath{{tikzfigs/}{figs/}}
\usepackage{subcaption}
\usepackage{array}
\usepackage{color}
\usepackage{hhline}
\usepackage{multirow}
\usepackage{enumitem}
\usepackage{bm}
\usepackage{ulem}
\usepackage{booktabs}
\usepackage{algorithm}
\usepackage{algpseudocode}
\usepackage{lipsum}

\usepackage[version=4]{mhchem}
\usepackage{siunitx}
\usepackage{longtable,tabularx}
\makeatletter
\newcommand*\bigcdot{\mathpalette\bigcdot@{.4}}
\newcommand*\bigcdot@[2]{\mathbin{\vcenter{\hbox{\scalebox{#2}{$\m@th#1\bullet$}}}}}
\makeatother

\newcommand{\sref}[1]{\S\ref{#1}}
\newcommand{\oursec}[1]{\vspace{0.1cm}\noindent{\normalsize\textit{#1}:}\,}

\newtheorem{remark}{Remark}

\newcommand{\boxing}[6]{
	\begin{figure*}[#1]
		\noindent\makebox[\textwidth][c] {
			\fbox{
				\begin{minipage}{#4}
					\vspace{0.2cm}
					\begin{#2} {\bf \textit{#5}} \label{#3} \end{#2}
					#6
				\end{minipage}
			}
		}
	\end{figure*}
}

\newcommand{\bvec}[1]{#1}%{\bm{#1}}

\newcommand{\definedas}{\coloneqq}
\newcommand{\dg}{^{\circ}}

\newcommand{\skewsymbig}{\Omega}
\newcommand{\OMEGA}[1]{\Omega\big(#1\big)}
\newcommand{\diag}[1]{\mathop{\mathrm{diag}}\left(#1\right)}
\newcommand{\fcn}[3]{#1:#2 \rightarrow #3}
\newcommand{\pder}[2]{\frac{\partial #1}{\partial #2}}
\newcommand{\dd}[2]{\frac{d#1}{d#2}}
\newcommand{\dotprod}[2]{#1 \bigcdot #2}
\newcommand{\crossprod}[2]{#1 \times #2}
\newcommand{\ith}[2]{#1^{\textit{#2}}}

\newcommand{\intee}[2]{[#1\,,#2]}
\newcommand{\intei}[2]{[#1\,,#2)}
\newcommand{\intie}[2]{(#1\,,#2]}
\newcommand{\intii}[2]{(#1\,,#2)}

\newcommand{\real}{\mathbb{R}}
\newcommand{\sthree}{\mathcal{S}^3}
\newcommand{\sothree}{\text{SO}(3)}
\newcommand{\spd}[1]{\mathbb{S}^{#1}_{++}}

\newcommand{\onenorm}[1]{\left\|#1\right\|_1}
\newcommand{\twonorm}[1]{\left\|#1\right\|_2}

\newcommand{\refframe}{\mathcal{F}}
\newcommand{\inertial}{\mathcal{I}}
\newcommand{\body}{\mathcal{B}}
\newcommand{\iframe}{\refframe_\inertial}
\newcommand{\bframe}{\refframe_\body}
\newcommand{\ex}{\bvec{e}_{1}}
\newcommand{\ey}{\bvec{e}_{2}}
\newcommand{\ez}{\bvec{e}_{3}}
\newcommand{\qidentity}{\bvec{q}_{\textit{id}}}

\newcommand{\subin}{\textit{in}}
\newcommand{\subig}{\textit{ig}}
\newcommand{\subf}{\textit{f}}
\newcommand{\tin}{t_{\subin}} % initial time
\newcommand{\tig}{t_{\subig}} % ignition time
\newcommand{\tf}{t_{\subf}} % final time
\newcommand{\tc}{t_{\textit{c}}} % coast time
\newcommand{\tb}{t_{\textit{b}}} % burn time
\newcommand{\tcmax}{t_{\textit{c},\textit{max}}} % longest allowable coast time

\newcommand{\m}{m}
\newcommand{\pos}{\bvec{r}}
\newcommand{\vel}{\bvec{v}}
\newcommand{\thrust}{\bvec{T}}
\newcommand{\aero}{\bvec{A}}
\newcommand{\rI}{\bvec{r}_\inertial}
\newcommand{\vI}{\bvec{v}_\inertial}
\newcommand{\vB}{\bvec{v}_\body}
\newcommand{\qIB}{\bvec{q}_{\body\leftarrow\inertial}}
\newcommand{\omegaB}{\bvec{\omega}_\body}
\newcommand{\TI}{\bvec{T}_\inertial}
\newcommand{\TB}{\bvec{T}_\body}
\newcommand{\aeroI}{\bvec{A}_\inertial}
\newcommand{\aeroB}{\bvec{A}_\body}
\newcommand{\xx}{\bvec{x}} % state vector
\newcommand{\uu}{\bvec{u}} % control vector

\newcommand{\mdot}{\dot{m}}
\newcommand{\rIdot}{\dot{\bvec{r}}_\inertial}
\newcommand{\vIdot}{\dot{\bvec{v}}_\inertial}
\newcommand{\qIBdot}{\dot{\bvec{q}}_{\body\leftarrow\inertial}}
\newcommand{\omegaBdot}{\dot{\bvec{\omega}}_\body}

\newcommand{\rIin}{\bvec{r}_{\inertial,\subin}}
\newcommand{\vIin}{\bvec{v}_{\inertial,\subin}}

\newcommand{\mig}{\m_{\subig}}
\newcommand{\rIig}{\bvec{r}_{\inertial,\subig}}
\newcommand{\vIig}{\bvec{v}_{\inertial,\subig}}

\newcommand{\omegaBig}{\bvec{0}}

\newcommand{\prig}{p_{\bvec{r},\subig}}
\newcommand{\pvig}{p_{\bvec{v},\subig}}

\newcommand{\rIf}{\bvec{0}}
\newcommand{\vdes}{v_{d}}
\newcommand{\vIf}{-\vdes\bvec{e}_{1}}
\newcommand{\qIBf}{\qidentity}
\newcommand{\omegaBf}{\bvec{0}}

\newcommand{\Isp}{I_{\textit{sp}}}
\newcommand{\Anoz}{A_{\textit{noz}}}
\newcommand{\mdotalpha}{\alpha_{\dot{\m}}}
\newcommand{\mdotbeta}{\beta_{\dot{\m}}}
\newcommand{\rTB}{\bvec{r}_{T,\body}}
\newcommand{\rCPB}{\bvec{r}_{\textit{cp},\body}}
\newcommand{\inertia}{J_\body}
\newcommand{\Ca}{C_{A}}

\newcommand{\Sa}{S_{A}}
\newcommand{\setaero}{\mathcal{A}(V)}

\newcommand{\tiltmax}{\theta_{\textit{max}}}
\newcommand{\omegamax}{\omega_{\textit{max}}}
\newcommand{\glideslope}{\gamma_{\textit{gs}}}
\newcommand{\mdry}{m_{\textit{dry}}}
\newcommand{\gimbalmax}{\delta_{\textit{max}}}
\newcommand{\Tmin}{T_{\textit{min}}}
\newcommand{\Tmax}{T_{\textit{max}}}

\newcommand{\stcnx}{{n_z}}

\newcommand{\stcg}{g} % STC trigger function
\newcommand{\stcf}{c} % STC constraint function
\newcommand{\stch}{h} % STC mixed function
\newcommand{\stcs}{\sigma} % STC slack variable
\newcommand{\stcshat}{\hat{\stcs}} % STC LCP slack variable solution
\newcommand{\stcx}{\bvec{z}} % generalized optimization vector
\newcommand{\qalpha}{$q$-$\aoa$}
\newcommand{\aoa}{\alpha}
\newcommand{\Vaoa}{V_{\aoa}}
\newcommand{\aoamax}{\aoa_{\textit{max}}}
\newcommand{\Rfov}{R_{\textit{fov}}}
\newcommand{\losmax}{\lambda_{\textit{max}}}
\newcommand{\gaoa}[1]{\stcg_{\aoa}\big(#1\big)}
\newcommand{\faoa}[2]{\stcf_{\aoa}\big(#1,#2\big)}
\newcommand{\haoa}[2]{\stch_{\aoa}\big(#1,#2\big)}

\newcommand{\gstd}{g_0}
\newcommand{\gI}{\bvec{g}_\inertial}
\newcommand{\density}{\rho}
\newcommand{\Pamb}{P_{\textit{amb}}}
\newcommand{\cIB}{C_{\body\leftarrow\inertial}}
\newcommand{\cBI}{C_{\inertial\leftarrow\body}}
\newcommand{\FI}{\bvec{F}_\inertial}
\newcommand{\MB}{\bvec{M}_\body}
\newcommand{\Hgs}{H_{\gamma}}
\newcommand{\Htilt}{H_{\theta}}

\newcommand{\ns}{{n_x}} % dimension of state vector
\newcommand{\nc}{{n_u}} % dimension of control vector
\newcommand{\nz}{{n_z}} % dimension of state-control-parameter vector
\newcommand{\np}{{n_{\textit{cvx}}}} % number of convex constraints
\newcommand{\nq}{{n_{\textit{ncvx}}}} % number of non-convex constraints
\newcommand{\nr}{{n_{\textit{stc}}}} % number of cSTC constraints
\newcommand{\dilation}{s}
\newcommand{\zz}{\bvec{z}} % stacked state-control-parameter vector
\newcommand{\rr}{\bvec{w}} % linearization residual vector
\newcommand{\tco}{\bar{t}_{c}} % linearization initial time
\newcommand{\sso}{\bar{\dilation}} % linearization dilation
\newcommand{\xxo}{\bar{\xx}} % reference state vector
\newcommand{\uuo}{\bar{\uu}} % reference control vector
\newcommand{\zzo}{\bar{\zz}} % reference stacked state-control-parameter vector
\newcommand{\fcnJ}{J} % cost function
\newcommand{\fcnf}{f} % nonlinear dynamics
\newcommand{\fcnsf}{F}%\tilde{f}} % temporally normalized nonlinear dynamics
\newcommand{\fcncons}{h} % equality constraints
\newcommand{\Jxu}[2]{\fcnJ\big(#1\,,#2\big)}
\newcommand{\fxu}[2]{\fcnf\big(#1\,,#2\big)}
\newcommand{\fxus}[3]{\fcnsf\big(#1\,,#2\,,#3\big)}
\newcommand{\setcvx}{\mathcal{I}_{\textit{cvx}}}
\newcommand{\setncvx}{\mathcal{I}_{\textit{ncvx}}}
\newcommand{\setcstc}{\mathcal{I}_{\textit{stc}}}
\newcommand{\tnorm}{\tau}
\newcommand{\Ac}{A}
\newcommand{\Bc}{B}
\newcommand{\Sc}{\bvec{S}}
\newcommand{\dfsdx}{\displaystyle{\pder{\fcnsf}{\xx}\Bigr|_{\zzo(\tnorm)}}}
\newcommand{\dfsdu}{\displaystyle{\pder{\fcnsf}{\uu}\Bigr|_{\zzo(\tnorm)}}}
\newcommand{\dfsds}{\displaystyle{\pder{\fcnsf}{\dilation}\Bigr|_{\zzo(\tnorm)}}}
\newcommand{\hi}[2]{\fcncons_{#2}\big(#1\big)}
\newcommand{\dhidz}{\displaystyle{\pder{\fcncons_i}{\zz}\Bigr|_{\zzo(\tnorm)}}}

\newcommand{\KK}{K}
\newcommand{\setK}{\mathcal{K}}
\newcommand{\setKm}{\bar{\mathcal{K}}}
\newcommand{\tauk}{\tau_k}
\newcommand{\taukp}{\tau_{k+1}}
\newcommand{\xxk}{\xx_k}
\newcommand{\uuk}{\uu_k}
\newcommand{\zzk}{\zz_k}
\newcommand{\xxko}{\xxo_k}
\newcommand{\uuko}{\uuo_k}
\newcommand{\zzko}{\zzo_k}
\newcommand{\xxkp}{\xx_{k+1}}
\newcommand{\uukp}{\uu_{k+1}}
\newcommand{\lambdalk}[1]{\lambda_{k}^{-}(#1)}
\newcommand{\lambdark}[1]{\lambda_{k}^{+}(#1)}
\newcommand{\stm}[2]{\Phi_{A}\big(#1\,,#2\big)}
\newcommand{\stminv}[2]{\Phi_{A}^{-1}\big(#1\,,#2\big)}
\newcommand{\Ad}{A_{k}}
\newcommand{\Bdm}{B^{-}_{k}}
\newcommand{\Bdp}{B^{+}_{k}}
\newcommand{\Sd}{\bvec{S}_{k}}
\newcommand{\rd}{\rr_k}

\newcommand{\nuk}{\bvec{\nu}_k}
\newcommand{\nubig}{\bvec{\nu}}
\newcommand{\wvc}{w_{\nu}}
\newcommand{\Wtr}{W_{\textit{tr}}}
\newcommand{\Jvc}[1]{J_{\textit{vc}}(#1)}
\newcommand{\Jtr}[2]{J_{\textit{tr}}(#1\,,#2)}
\newcommand{\epsvc}{\epsilon_{\textit{vc}}}
\newcommand{\epstr}{\epsilon_{\textit{tr}}}

\newcommand{\mki}{\m_{1}}
\newcommand{\rIki}{\bvec{r}_{\inertial,1}}
\newcommand{\vIki}{\bvec{v}_{\inertial,1}}

\newcommand{\omegaBki}{\bvec{\omega}_{\body,1}}

\newcommand{\xxoig}{\xxo_{ig}}

\newcommand{\mk}{\m_{k}}
\newcommand{\rIk}{\bvec{r}_{\inertial,k}}

\newcommand{\qIBk}{\bvec{q}_{\body\leftarrow\inertial,k}}
\newcommand{\omegaBk}{\bvec{\omega}_{\body,k}}

\newcommand{\xxkf}{\xx_{K}}
\newcommand{\mkf}{\m_{K}}
\newcommand{\rIkf}{\bvec{r}_{\inertial,K}}
\newcommand{\vIkf}{\bvec{v}_{\inertial,K}}
\newcommand{\qIBkf}{\bvec{q}_{\body\leftarrow\inertial,K}}
\newcommand{\omegaBkf}{\bvec{\omega}_{\body,K}}

\newcommand{\xxof}{\xxo_{f}}

\newcommand{\ftlbk}{h_{\textit{tlb},k}}
\newcommand{\Htlbk}{H_{\textit{tlb},k}}

\newcommand{\haoak}{h_{\aoa,k}}
\newcommand{\Haoak}{H_{\aoa,k}}

\newcommand{\LU}{U_L}
\newcommand{\TU}{U_T}
\newcommand{\MU}{U_M}

\usepackage{rotating}
\usepackage{tikz}
\usetikzlibrary{patterns}
\usetikzlibrary{math}
\usetikzlibrary{scopes}
\usetikzlibrary{fadings}
\usetikzlibrary{arrows,calc,decorations.pathmorphing}
\usetikzlibrary{matrix,positioning}
\tikzset{>=latex}

\definecolor{beige}{RGB}{245,245,220}
\definecolor{darkred}{rgb}{0.90,0.00,0.00}%red}%!75!black!100}
\definecolor{darkgreen}{rgb}{0.00,0.45,0.00}

\tikzset{ 
table/.style={
  matrix of math nodes,
  row sep=-\pgflinewidth,
  column sep=-\pgflinewidth,
  nodes={rectangle,draw=white,text width=6em,align=left},
  text depth=0.25ex,
  text height=2ex,
  nodes in empty cells
  },
  title/.style={font=\large}
}

\definecolor{ucol}{RGB}{255,0,0}
\definecolor{gcol}{RGB}{0,120,0}
\definecolor{scol}{RGB}{63, 226, 45}
\definecolor{vcol}{RGB}{0,0,0}
\definecolor{acol}{RGB}{0, 128, 255}
\definecolor{dcol}{RGB}{204,102,0}
\definecolor{lcol}{RGB}{204,102,0}
\definecolor{bcol}{RGB}{0,0,0}
\definecolor{ocol}{RGB}{167,167,167}

\newcommand{\threeaxes}[8]{
	\tikzmath{
    	\Lxr=#3;\Lxl=#3;\Lyt=#4;\Lyb=#4;\Zt= #5;\Zb= #6;\Xang=#7;\Yang=#8;
    	\Xxr= cos(\Xang)*\Lxr; \Xyr=-sin(\Xang)*\Lxr;
    	\Xxl=-cos(\Xang)*\Lxl; \Xyl= sin(\Xang)*\Lxl;
        \Yxt= sin(\Yang)*\Lyt; \Yyt= cos(\Yang)*\Lyt;
        \Yxb=-sin(\Yang)*\Lyb; \Yyb=-cos(\Yang)*\Lyb;
        \zzz=0;
    }
    \begin{scope}[shift={(#1,#2)},rotate=0]
        \ifx\Zt\zzz\else  \draw[black,->] (0,0) -- +(0, \Zt);    \fi
        \ifx\Zb\zzz\else  \draw[black,->] (0,0) -- +(0,-\Zb);    \fi
        \ifx\Lxr\zzz\else \draw[black,->] (0,0) -- +(\Xxr,\Xyr); \fi
        \ifx\Lxl\zzz\else \draw[black,->] (0,0) -- +(\Xxl,\Xyl); \fi
        \ifx\Lyt\zzz\else \draw[black,->] (0,0) -- +(\Yxt,\Yyt); \fi
        \ifx\Lyb\zzz\else \draw[black,->] (0,0) -- +(\Yxb,\Yyb); \fi
	\end{scope}
}
\newcommand{\threeaxeslabelx}[9]{
	\tikzmath{
    	\Lxr=#3;\Lxl=#3;\Lyt=#4;\Lyb=#4;\Zt= #5;\Zb= #6;\Xang=#7;\Yang=#8;
    	\Xxr= cos(\Xang)*\Lxr; \Xyr=-sin(\Xang)*\Lxr;
    	\Xxl=-cos(\Xang)*\Lxl; \Xyl= sin(\Xang)*\Lxl;
        \Yxt= sin(\Yang)*\Lyt; \Yyt= cos(\Yang)*\Lyt;
        \Yxb=-sin(\Yang)*\Lyb; \Yyb=-cos(\Yang)*\Lyb;
        \zzz=0;
    }
    \begin{scope}[shift={(#1,#2)},rotate=0]
        \ifx\Lxr\zzz\else \draw (\Xxr,\Xyr) node[anchor=south west] {#9}; \fi
	\end{scope}
}
\newcommand{\threeaxeslabely}[9]{
	\tikzmath{
    	\Lxr=#3;\Lxl=#3;\Lyt=#4;\Lyb=#4;\Zt= #5;\Zb= #6;\Xang=#7;\Yang=#8;
    	\Xxr= cos(\Xang)*\Lxr; \Xyr=-sin(\Xang)*\Lxr;
    	\Xxl=-cos(\Xang)*\Lxl; \Xyl= sin(\Xang)*\Lxl;
        \Yxt= sin(\Yang)*\Lyt; \Yyt= cos(\Yang)*\Lyt;
        \Yxb=-sin(\Yang)*\Lyb; \Yyb=-cos(\Yang)*\Lyb;
        \zzz=0;
    }
    \begin{scope}[shift={(#1,#2)},rotate=0]
        \ifx\Lyt\zzz\else \draw (\Yxt,\Yyt) node[anchor=south west] {#9}; \fi
	\end{scope}
}
\newcommand{\threeaxeslabelz}[9]{
	\tikzmath{
    	\Lxr=#3;\Lxl=#3;\Lyt=#4;\Lyb=#4;\Zt= #5;\Zb= #6;\Xang=#7;\Yang=#8;
    	\Xxr= cos(\Xang)*\Lxr; \Xyr=-sin(\Xang)*\Lxr;
    	\Xxl=-cos(\Xang)*\Lxl; \Xyl= sin(\Xang)*\Lxl;
        \Yxt= sin(\Yang)*\Lyt; \Yyt= cos(\Yang)*\Lyt;
        \Yxb=-sin(\Yang)*\Lyb; \Yyb=-cos(\Yang)*\Lyb;
        \zzz=0;
    }
    \begin{scope}[shift={(#1,#2)},rotate=0]
        \ifx\Zt\zzz\else  \draw (0,\Zt)     node[anchor=west] {$\;$#9}; \fi
	\end{scope}
}

\newcommand{\isoaxes}[9]{
    \tikzmath{
        \rot=#3;
        \len=#4;
        \ddd=#5;
        \pX=  0; \Xx=cos(\pX)*\len; \Xy=sin(\pX))*\len;
        \pY=120; \Yx=cos(\pY)*\len; \Yy=sin(\pY))*\len;
        \pZ=240; \Zx=cos(\pZ)*\len; \Zy=sin(\pZ))*\len;
        \Xxx=cos(\ddd)*\Xx-sin(\ddd)*\Xy; \Xxy=sin(\ddd)*\Xx+cos(\ddd)*\Xy;
        \Yxx=cos(\ddd)*\Yx-sin(\ddd)*\Yy; \Yxy=sin(\ddd)*\Yx+cos(\ddd)*\Yy;
        \Zxx=cos(\ddd)*\Zx-sin(\ddd)*\Zy; \Zxy=sin(\ddd)*\Zx+cos(\ddd)*\Zy;
    }
    \begin{scope}[shift={(#1,#2)},rotate=\rot]
		\filldraw[black] (0,0) circle (2pt);
        \draw[black,thick,->] (0,0) -- +(\Xx,\Xy);
        \draw[black,thick,->] (0,0) -- +(\Yx,\Yy);
        \draw[black,thick,->] (0,0) -- +(\Zx,\Zy);
        \draw (0.1,0.6)   node[rotate=0,anchor=center] {#6};
        \draw (\Xxx,\Xxy) node[rotate=0,anchor=center] {#7};
        \draw (\Yxx,\Yxy) node[rotate=0,anchor=center] {#8};
        \draw (\Zxx,\Zxy) node[rotate=0,anchor=center] {#9};
    \end{scope}
}

\newcommand{\centerofmass}[1]{%
    \tikz[radius=0.4em,scale=#1] {%
        \fill (0,0) -- ++(0.4em,0) arc [start angle=0,end angle=90] -- ++(0,-0.8em) arc [start angle=270, end angle=180];%
        \draw (0,0) circle;%
    }%
}

\newcommand{\pane}[7]{
	\tikzmath{
    	\rot=#3;
    	\width=#4;
        \height=#5;
        \corner=#6;
    	\px1= 0.5*\width-\corner; \py1=-0.5*\height;
        \px2= 0.5*\width;         \py2=-0.5*\height+\corner;
        \px3= 0.5*\width;         \py3= 0.5*\height-\corner;
        \px4= 0.5*\width-\corner; \py4= 0.5*\height;
        \px5=-0.5*\width+\corner; \py5= 0.5*\height;
        \px6=-0.5*\width;         \py6= 0.5*\height-\corner;
        \px7=-0.5*\width;         \py7=-0.5*\height+\corner;
        \px8=-0.5*\width+\corner; \py8=-0.5*\height;
    }
	\begin{scope}[shift={(#1,#2)},rotate=\rot]
    	\filldraw[{#7}]
        (\px1,\py1) to[out=    0,in=  -90] (\px2,\py2) --
        (\px3,\py3) to[out=   90,in=    0] (\px4,\py4) --
        (\px5,\py5) to[out= -180,in=   90] (\px6,\py6) --
        (\px7,\py7) to[out=  -90,in= -180] (\px8,\py8) -- cycle;
    \end{scope}
}
\newcommand{\cpane}[7]{
	\tikzmath{
		\pLTx=#1;
		\pRBx=#2;
		\pLTy=#3;
		\pRBy=#4;
	}
	\pane{0.5*\pLTx+0.5*\pRBx}{0.5*\pLTy+0.5*\pRBy}{#5}{\pRBx-\pLTx}{\pLTy-\pRBy}{#6}{#7}
}

\newcommand{\coneback}[7]{
	\tikzmath{\rot=#3;
              \length=#4; \radius=\length*tan(0.5*#5); \depth=#6;
              \sx =  cos(\rot)*#1 + sin(\rot)*#2;
              \sy = -sin(\rot)*#1 + cos(\rot)*#2;
    }
    \begin{scope}[shift={(\sx,\sy)},transform canvas={rotate=\rot}]
	    \draw[{#7}] (\radius,-\length) arc(360:180: {\radius} and {-\depth});
    \end{scope}
}

\newcommand{\cone}[7]{
	\tikzmath{\rot=#3;
              \length=#4; \radius=\length*tan(0.5*#5); \depth=#6;
              \sx =  cos(\rot)*#1 + sin(\rot)*#2;
              \sy = -sin(\rot)*#1 + cos(\rot)*#2;
	}
    \begin{scope}[shift={(\sx,\sy)},transform canvas={rotate=\rot}]
    	\fill[{#7}] (0,0) -- (\radius,-\length) arc(360:180: {\radius} and {\depth}) -- cycle;
    	\draw[color=black!100] (0,0) -- (\radius,-\length) arc(360:180: {\radius} and {\depth}) -- cycle;
    \end{scope}
}

\newcommand{\rocket}[6]{
	\tikzmath{\rot=#3;
    		  \length=#4;
              \throttle=(\length/0.8)*0.6*#5;
              \gimbalangle=#6;
    		  \LL = \length; \RR = \LL/8;
              \HH = \LL/5;   \WW = \LL/16;
              \rrrr = \LL/16;  \ww = \LL/16;
              \HG = \LL/8;   \RG = \LL/11; \WG = \LL/16;
              \sx =  cos(\rot)*#1 + sin(\rot)*#2;
              \sy = -sin(\rot)*#1 + cos(\rot)*#2;
    }
    \begin{scope}[shift={(\sx,\sy)},transform canvas={rotate=\rot}]
		\begin{scope}[shift={(0,-0.5*\LL)},rotate=\gimbalangle]
          	\fill[fill=orange!100,shading=axis,shading angle=90,left color=orange!100,right color=orange!25]
            	(-\RG,-\HG) -- (0,-\HG-\throttle) -- (\RG,-\HG) -- cycle;
        	\filldraw[color=black!100,fill=black!10,shading=axis,shading angle=90,left color=black!30,right color=black!0]
        	(0,0) -- (\RG,-\HG) arc(360:180: {\RG} and {\WG}) -- cycle;
		\end{scope}
        
    	\filldraw[color=black!100,fill=black!10,shading=ball,shading angle=90,left color=black!30,right color=black!0]
        	(-\RR,-0.5*\LL) arc(180:360: {\RR} and {\WW}) -- (\RR,0.5*\LL) -- (\RR,0.5*\LL) arc(360:180: {\RR} and {\WW}) -- (-\RR,-0.5*\LL) -- cycle;
        
    	\filldraw[color=black!100,fill=black!10,shading=axis,shading angle=90,left color=black!30,right color=black!0]
        	(\RR,0.5*\LL) arc(360:180: {\RR} and {\WW}) -- (-\rrrr,0.5*\LL+\HH) arc(-180:0: {\rrrr} and {-\ww}) -- cycle;
    \end{scope}
}

\title{Successive Convexification for Real-Time 6-DoF Powered Descent Guidance with State-Triggered Constraints}

\author{Michael Szmuk\footnote{PhD Candidate, UW Aero. and Astro., mszmuk@uw.edu, AIAA Member.}, Taylor P. Reynolds\footnote{PhD Candidate, UW Aero. and Astro., tpr6@uw.edu, AIAA Member.}, and Beh\c{c}et A\c {c}{\i}kme\c{s}e\footnote{Associate Professor, UW Aero. and Astro., behcet@uw.edu, AIAA Member.}}
\affil{University of Washington, Seattle, WA, 98195-2400}

\begin{document}
\maketitle

%%%%%%%%%%%%%%%%%%%%%%%%%%%%%%%%%%%%%%%%%%%%%%%%%%%%%%%%%%%%%%%%%%%%%%%%%%%%%%%%%%%%%%%
%%%%%%%%%%%%%%%%%%%%%%%%%%%%%%%%%%%%%%%%%%%%%%%%%%%%%%%%%%%%%%%%%%%%%%%%%%%%%%%%%%%%%%%
%%%%%%%%%%%%%%%%%%%%%%%%%%%%%%%%%%%%%%%%%%%%%%%%%%%%%%%%%%%%%%%%%%%%%%%%%%%%%%%%%%%%%%%

\begin{abstract}

In this paper, we present a real-time successive convexification algorithm for a generalized free-final-time 6-degree-of-freedom powered descent guidance problem. We build on our previous work by introducing the following contributions: (i) a free-ignition-time modification that allows the algorithm to determine the optimal engine ignition time, (ii) a tractable aerodynamics formulation that models both lift and drag, and (iii) a continuous state-triggered constraint formulation that emulates conditionally enforced constraints. In particular, contribution (iii) effectively allows constraints to be enabled or disabled by \textit{if}-statements conditioned on the solution variables of the parent continuous optimization problem. To the best of our knowledge, this represents a novel formulation in the optimal control literature, and enables a number of interesting applications, including velocity-triggered angle of attack constraints and range-triggered line of sight constraints. Our algorithm converts the resulting generalized powered descent guidance problem from a non-convex free-final-time optimal control problem into a sequence of tractable convex second-order cone programming subproblems. With the aid of virtual control and trust region modifications, these subproblems are solved in succession until convergence is attained. Simulations using a third-party solver demonstrate the real-time capabilities of the proposed algorithm, with a maximum execution time of less than~$0.7$ seconds over a multitude of problem feature combinations.

\end{abstract}

%%%%%%%%%%%%%%%%%%%%%%%%%%%%%%%%%%%%%%%%%%%%%%%%%%%%%%%%%%%%%%%%%%%%%%%%%%%%%%%%%%%%%%%
%%%%%%%%%%%%%%%%%%%%%%%%%%%%%%%%%%%%%%%%%%%%%%%%%%%%%%%%%%%%%%%%%%%%%%%%%%%%%%%%%%%%%%%
%%%%%%%%%%%%%%%%%%%%%%%%%%%%%%%%%%%%%%%%%%%%%%%%%%%%%%%%%%%%%%%%%%%%%%%%%%%%%%%%%%%%%%%

\section{Introduction} \label{sec:intro}

\lettrine{T}{his} paper presents a real-time guidance algorithm that solves a generalized free-final-time 6-degree-of-freedom~(DoF) powered descent guidance problem with free ignition time, aerodynamic effects, and conditionally enforced constraints. Real-time optimal guidance algorithms are an enabling technology for future manned and unmanned planetary missions that require autonomous precision landing capabilities. Such algorithms allow for the explicit inclusion of operational and mission constraints, and are able to compute trajectories that are (locally) optimal with respect to key metrics, such as propellant consumption, burn time, or miss distance. Consequently, these algorithms enhance a lander's ability to recover from a wider range of dispersions encountered during the entry, descent, and landing phase, and to react to obstructions on the surface that become apparent only as the vehicle approaches the landing site.

The constraint satisfaction afforded by these algorithms allows designers to select more ambitious and scientifically interesting landing sites, enhances the lander's ability to handle uncertainties, and ultimately increases the likelihood of mission success. Moreover, the optimality may be leveraged to improve the scientific return of planetary missions by reducing the propellant mass fraction~\cite{Scharf15}. The relevance of real-time optimal guidance algorithms been demonstrated in the (now routine) landings of orbital-class vertical-takeoff-vertical-landing  reusable launch vehicles~\cite{larsNAE}.

Solving the 6-DoF powered descent guidance problem in real-time is challenging for several reasons. First, the problem consists of nonlinear dynamics and non-convex state and control constraints, and does not yet have an analytical solution. Second, since the problem must be solved using numerical methods, the validity of the solution is depends heavily on the discretization scheme. Third, the non-convex nature of the problem makes it difficult to select a suitable reference trajectory to initialize an iterative solution process. The successive convexification methodology used in this paper addresses these issues, and is able to solve the problem quickly and reliably over a wide range of conditions.

\subsection{Related Work}

The powered descent guidance literature can be separated into works that consider 3-DoF translation dynamics, and ones that consider more general 6-DoF rigid body dynamics. Work on 3-DoF powered descent guidance began many years ago during the Apollo program, with several authors approaching the problem using optimal control theory~\cite{Meditch1964,Lawden1963} and the calculus of variations~\cite{Marec1979}. These early works noted that the solution to the fuel-optimal powered descent guidance problem exhibited bang-bang behavior, where the thrust assumed either the minimum or maximum allowable value everywhere along the optimal trajectory. While these results offered important insights into the powered descent guidance problem, they were not incorporated into the Apollo flight code, since the polynomial-based guidance methods used to perform the landings were deemed sufficiently close to optimal, and were far simpler to design~\cite{klumpp}. In the years following Apollo, researchers continued searching for analytical solutions to the 3-DoF landing problem~\cite{Kornhauser1972,Breakwell1975,Azimov1996,dsouza97} that would be conducive to on-board implementation for future missions.

Unmanned missions to the Martian surface in the early the 2000s renewed interest in using direct methods to compute solutions to the powered descent guidance problem. In 2005, Topcu, Casoliva, and Mease presented results for the 3-DoF problem that resembled a modern take on~\cite{Lawden1963}, adding numerical simulations to reinforce and demonstrate theoretical results~\cite{Topcu2005,topcu_pdg}. Around the same time, A\c{c}{\i}kme\c{s}e and Ploen published work on a convex programming approach to the problem~\cite{Acikmese2005,ploenaiaa06,behcetjgcd07}. Using Pontryagin's maximum principle, they showed that the non-convex thrust magnitude lower bound constraint could be losslessly convexified by introducing a relaxation that rendered the optimal solution of the (now convex) relaxed problem identical to that of the original non-convex problem. Hence, the difficult task of solving the non-convex 3-DoF powered descent guidance problem was converted into the far simpler, but equivalent, task of solving a convex second-order cone programming problem. Subsequently, lossless convexification was extended to include non-convex pointing constraints, and to encompass more general optimal control problems~\cite{behcetaut11,harris_acc,lars_12scl,behcet_aut11,larssys12,pointing2013,matt_aut1}. This methodology was demonstrated in a sequence of flight experiments in the early 2010s~\cite{dueri2016customized,Scharf14}.

More recently, researchers have considered the 6-DoF powered descent guidance problem. In 2012, Lee and Mesbahi formulated the powered descent guidance problem using dual quaternions~\cite{Lee2012}. This work revealed that certain coupled rotational-translational constraints are convex when using this parameterization~\cite{Lee2015}. In~\cite{Lee2017}, a piecewise-affine approximation of the dynamics was used to design a model predictive controller that generated thrust and torque commands. However, since the accuracy of the equations of motion relied heavily on the temporal resolution of the approximation, the approach resulted in either short prediction horizons or prohibitively large optimization problems.

Work on powered descent guidance and atmospheric entry problems turned to Sequential Convex Programming (SCP) methods in order to handle more general non-convexities~\cite{jordithesis,liu2014solving,liu2015entry,wang2016constrained}. These methods solve non-convex problems by iteratively solving a sequence of local convex approximations obtained via linearization. The first order approach used in SCP methods guarantees that the approximations are convex, and lies in stark contrast to the second order approach used in Sequential Quadratic Programming (SQP) methods~\cite{Boggs1995}, which may expend significant computational effort to ensure the convexity of each approximation.

In this paper, we solve a generalized free-final-time 6-DoF powered descent guidance problem using the successive convexification framework. Successive convexification can be classified as an SCP method that uses virtual control and trust region modifications to facilitate convergence. The algorithms detailed in~\cite{SCvx_cdc16,mao2016successive,SCvx_2017arXiv} use an exact penalty method in conjunction with hard trust regions, whereas those used in this paper and in~\cite{szmuk2016successive,szmuk2017successive,szmuk2018successive} employ soft trust regions.

\subsection{Statement of Contribution}

This paper presents three primary contributions: (i) a free-ignition-time modification that allows the algorithm to determine the optimal engine ignition time, (ii) a tractable aerodynamics formulation that models both lift and drag, and (iii) a continuous state-triggered constraint formulation that emulates conditionally enforced constraints. In particular, contribution (iii) effectively allows constraints to be enabled or disabled by \textit{if}-statements in a continuous optimization framework.  

To the best of our knowledge, state-triggered constraints bear the most resemblance to two existing approaches: mixed-integer programming, and complementarity constraints. The former approach implements discrete decisions explicitly using integer variables, whereas the latter formulates such decisions implicitly using continuous variables. Despite the existence of efficient branch and bound techniques, mixed-integer formulations can suffer from poor complexity~\cite{Richards2015}, and are not conducive to solving the generalized powered descent guidance problem in real-time. Like state-triggered constraints, complementarity constraints~\cite{Cottle1992,Heemels2000} permit completely continuous formulations that can be more efficient than mixed-integer approaches (see Section~2.3 in~\cite{biegler2014}). However, complementarity constraints represent bi-directional \textit{if-and-only-if}-statements, whereas state-triggered constraints represent more general uni-directional \textit{if}-statements. Therefore, we argue that the proposed continuous state-triggered constraints are a key building block that enable the formulation of a broader set of constraints.

The secondary contributions of this paper are an improved description of the successive convexification methodology used in~\cite{szmuk2018successive}, and timing results that demonstrate the real-time capabilities of the algorithm. This paper regards the powered descent guidance problem as a feedforward trajectory generation problem, and does not address the topic of feedback control or issues arising from trajectory re-computation.

\subsection{Outline}

In~\sref{sec:problem_statement}, we present the primary contributions of this paper in the context of a generalized powered descent guidance problem. In~\sref{sec:convex_formulation}, we detail the successive convexification procedure and algorithm. In~\sref{sec:numerical_results}, we present simulation results that highlight the paper's contributions, and timing results that demonstrate the real-time capabilities of the proposed algorithm. Lastly,~\sref{sec:conclusion} provides concluding remarks.

%%%%%%%%%%%%%%%%%%%%%%%%%%%%%%%%%%%%%%%%%%%%%%%%%%%%%%%%%%%%%%%%%%%%%%%%%%%%%%%%%%%%%%%
%%%%%%%%%%%%%%%%%%%%%%%%%%%%%%%%%%%%%%%%%%%%%%%%%%%%%%%%%%%%%%%%%%%%%%%%%%%%%%%%%%%%%%%
%%%%%%%%%%%%%%%%%%%%%%%%%%%%%%%%%%%%%%%%%%%%%%%%%%%%%%%%%%%%%%%%%%%%%%%%%%%%%%%%%%%%%%%

\section{Problem Statement} \label{sec:problem_statement}

In this section, we present a 6-DoF formulation for a generalized powered descent guidance problem in the presence of atmospheric effects. This section is organized as follows. In~\sref{sec2:assumptions} and~\sref{sec2:notation}, we introduce the assumptions and notation used in our formulation. In~\sref{sec2:6dof}, we present a baseline problem formulation. In~\sref{sec2:fit}, \sref{sec2:aero}, and~\sref{sec2:stc}, we discuss the primary contributions of the paper, namely the free-ignition-time modification, the aerodynamic models, and state-triggered constraints. Lastly, \sref{sec2:nonconvex} provides a statement of the non-convex generalized powered descent guidance problem.

%%%%%%%%%%%%%%%%%%%%%%%%%%%%%%%%%%%%%%%%%%%%%%%%%%%%%%%%%%%%%%%%%%%%%%%%%%%%%%%%%%%

\subsection{Assumptions} \label{sec2:assumptions}

Most powered descent maneuvers commence at speeds substantially below orbital velocities and within only a few kilometers of the landing site. Hence, we neglect the effects of planetary rotation and assume a uniform gravitational field. We assume that the vehicle is equipped with a single rocket engine that can be gimbaled symmetrically about two axes up to a maximum gimbal angle, but stress that other thruster configurations can be readily accommodated. Further, we assume that the engine can be throttled between fixed minimum and maximum thrusts, and that once the engine is ignited it remains on until the terminal condition is reached. To tailor our treatment to applications with non-negligible atmospheric effects, we assume that the ambient atmospheric density and pressure are constant, that the aerodynamic forces are governed by the simplified models detailed in~\sref{sec2:aero}, and that the center-of-pressure is fixed with respect to a body-fixed reference frame. Further, we neglect the effects of winds, noting that constant uniform wind profiles can be readily incorporated into our formulation, and account for thrust reduction induced by atmospheric back-pressure by assuming the affine mass depletion dynamics in~\cite{szmuk2018successive}. Lastly, to make the problem tractable, we neglect higher order phenomena such as elastic structural modes and fuel slosh, and model the vehicle as a rigid body with a constant body-fixed center-of-mass and inertia.

%%%%%%%%%%%%%%%%%%%%%%%%%%%%%%%%%%%%%%%%%%%%%%%%%%%%%%%%%%%%%%%%%%%%%%%%%%%%%%%%%%%%

\subsection{Notation} \label{sec2:notation}

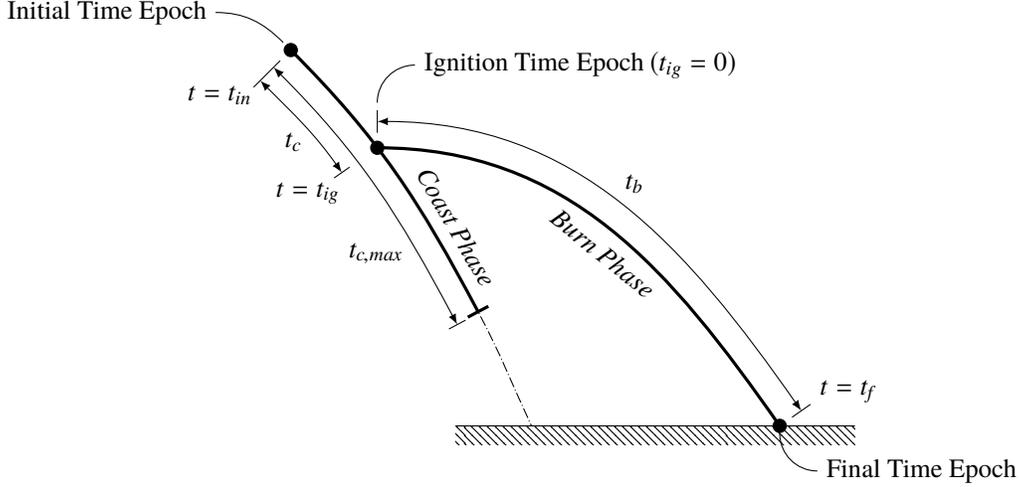
\begin{figure}[tb]
	\centering
    \begin{tikzpicture}
	\tikzset{>=latex}
    \tikzmath{\pAx=-0.5;  \pAy= 1.0;
		      \pBx= 2.0;  \pBy=-2.5;
			  \pCx= 0.65; \pCy=-0.3;
			  \pDx= 6.0;  \pDy=-4.0;
			  \pEx= 2.7;  \pEy= \pDy;
			  \dx=1.00;
			  \len1=0.2;
		      \len2=0.5+\len1;
			  \len3=0.75*\len1+0.25*\len2;
			  \len4=0.25*\len1+0.75*\len2;
			  \len5=0.3+\len1;
			  \len6=0.50*\len1+0.50*\len5;
			  \angA=45; \pAAx1=\pAx-\len1*sin(\angA); \pAAy1=\pAy-\len1*cos(\angA);
	                    \pAAx2=\pAx-\len2*sin(\angA); \pAAy2=\pAy-\len2*cos(\angA);
					    \pAAx3=\pAx-\len3*sin(\angA); \pAAy3=\pAy-\len3*cos(\angA);
				        \pAAx4=\pAx-\len4*sin(\angA); \pAAy4=\pAy-\len4*cos(\angA);
		      \angB=62; \pBBx1=\pBx-\len1*sin(\angB); \pBBy1=\pBy-\len1*cos(\angB);
					    \pBBx2=\pBx-\len2*sin(\angB); \pBBy2=\pBy-\len2*cos(\angB);
					    \pBBx3=\pBx-\len3*sin(\angB); \pBBy3=\pBy-\len3*cos(\angB);
					    \pBBx4=\pBx-\len4*sin(\angB); \pBBy4=\pBy-\len4*cos(\angB);
			  \angC=55; \pCCx1=\pCx-\len1*sin(\angC); \pCCy1=\pCy-\len1*cos(\angC);
			            \pCCx2=\pCx-\len2*sin(\angC); \pCCy2=\pCy-\len2*cos(\angC);
					    \pCCx3=\pCx-\len3*sin(\angC); \pCCy3=\pCy-\len3*cos(\angC);
					    \pCCx4=\pCx-\len4*sin(\angC); \pCCy4=\pCy-\len4*cos(\angC);
			  \angCC=0; \pCCCx1=\pCx+\len1*sin(\angCC); \pCCCy1=\pCy+\len1*cos(\angCC);
					    \pCCCx2=\pCx+\len5*sin(\angCC); \pCCCy2=\pCy+\len5*cos(\angCC);
					    \pCCCx3=\pCx+\len6*sin(\angCC); \pCCCy3=\pCy+\len6*cos(\angCC);
			 \angDD=55; \pDDDx1=\pDx+\len1*sin(\angDD); \pDDDy1=\pDy+\len1*cos(\angDD);
			            \pDDDx2=\pDx+\len5*sin(\angDD); \pDDDy2=\pDy+\len5*cos(\angDD);
			            \pDDDx3=\pDx+\len6*sin(\angDD); \pDDDy3=\pDy+\len6*cos(\angDD);
		     \pFx=0.5*\pAAx4+0.5*\pCCx4; \pFy=0.5*\pAAy4+0.5*\pCCy4;
			 \pGx=0.5*\pBBx4+0.5*\pCCx4; \pGy=0.5*\pBBy4+0.5*\pCCy4+0.2;
			 \pHx=0.5*\pCCCx3+0.5*\pDDDx3+0.5; \pHy=0.5*\pCCCy3+0.5*\pDDDy3+0.75;
		}
    
	    \definecolor{burncolor}{RGB}{0,0,0};
    	\definecolor{coastcolor}{RGB}{0,0,0};
%        \definecolor{burncolor}{RGB}{200,0,0};
%        \definecolor{coastcolor}{RGB}{0,120,0};
		\tikzset{dashdot/.style={dash pattern=on .4pt off 1pt on 4pt off 1pt}}
        
        % Coast Phase %
		\draw[coastcolor,very thick,-|] (\pAx,\pAy) to[out=-\angA,in=180-\angB] (\pBx,\pBy);
		\draw[black,<->] (\pAAx3,\pAAy3) to[out=-\angA,in=180-\angB] (\pBBx3,\pBBy3);
		\draw[black,<->] (\pAAx4,\pAAy4) to[out=-\angA,in=180-\angC] (\pCCx4,\pCCy4);
		\draw[black,dashdot] (\pBx,\pBy) to[out=-\angB,in=115] (\pEx,\pEy);
		\draw (\pAAx1,\pAAy1) -- (\pAAx2,\pAAy2);
		\draw (\pBBx1,\pBBy1) -- (0.5*\pBBx3+0.5*\pBBx4,0.5*\pBBy3+0.5*\pBBy4);
		\draw (0.5*\pCCx3+0.5*\pCCx4,0.5*\pCCy3+0.5*\pCCy4) -- (\pCCx2,\pCCy2);
		\draw (\pAAx2+0.1,\pAAy2-0.1) node[anchor=east] {$t=\tin$};
		\draw (\pCCx2+0.2,\pCCy2-0.2) node[anchor=east] {$t=\tig$};
		\draw (\pFx-0.1,\pFy-0.2) node[anchor=center] {$\tc$};
		\draw (\pGx-0.2,\pGy-0.2) node[anchor=center] {$\tcmax$};
		
		% Burn Phase %
		\draw[burncolor,very thick] (\pCx,\pCy) to[out=-0 ,in=180-\angDD] (\pDx,\pDy);
		\draw[black,<->] (\pCCCx3,\pCCCy3) to[out=-0 ,in=180-\angDD] (\pDDDx3,\pDDDy3);
		\draw (\pCCCx1,\pCCCy1) -- (\pCCCx2,\pCCCy2);
		\draw (\pDDDx1,\pDDDy1) -- (\pDDDx2,\pDDDy2);
		\draw (\pDDDx2,\pDDDy2-0.1) node[anchor=south west] {$t=\tf$};
		\draw (\pHx+0.1,\pHy+0.1) node[anchor=south] {$\tb$};
		
		% Ground
		\fill[pattern=north west lines, pattern color=black] (\pEx-\dx,\pEy) rectangle (\pDx+\dx,\pDy-0.25);
		\draw[black] (\pEx-\dx,\pEy) -- (\pDx+\dx,\pEy);
		
		% Points %
    	\filldraw[coastcolor] (\pAx,\pAy) circle (2.5pt);
		\filldraw[burncolor] (\pCx,\pCy) circle (2.5pt);
		\filldraw[burncolor] (\pDx,\pDy) circle (2.5pt);
		
		% Labels %
		\draw (\pAx-0.1,\pAy+0.1) to[out= 135,in=0  ] (\pAx-1.0,\pAy+0.5) node[anchor=east] {Initial Time Epoch};
		\draw (\pCCCx2,\pCCCy2+0.1) to[out= 90,in=-180] (\pCCCx2+0.5,\pCCCy2+0.6) node[anchor=west] {Ignition Time Epoch $(\tig = 0)$};
		\draw (\pDx,\pDy-0.1) to[out=-90,in=180] (\pDx+0.5,\pDy-0.6) node[anchor=west] {Final Time Epoch};
		\draw (\pCx+1.05,\pCy-1.05) node[anchor=center,rotate=-60,color=coastcolor] {\textit{Coast Phase}};
		\draw (\pCx+3.0,\pCy-1.4) node[anchor=center,rotate=-39,color=burncolor] {\textit{Burn Phase}};
\end{tikzpicture}
    \caption{Timeline of the generalized powered descent guidance problem. The problem begins at the initial time epoch~$\tin$ at a prescribed state and with the engine off. The vehicle is allowed to coast a maximum duration of~$\tcmax$ before engine ignition occurs at the ignition time epoch~$\tig$. The problem ends at the final time epoch~$\tf$, when the vehicle lands at the prescribed landing site.}
    \label{fig:traj_timeline}
\end{figure}

We begin by denoting time as~$t\in\real$, and define the \textit{initial time}~$\tin$ as the time at which the optimal control problem begins, the \textit{ignition time}~$\tig$ as the time at which the engine turns on, and the \textit{final time}~$\tf$ as the time at which the vehicle reaches the terminal condition. These epochs are defined such that~$\tin\leq\tig < \tf$, and their subscripts are used to denote problem parameters associated with the respective time epochs. Further, we define \textit{coast time} as~$\tc\definedas\tig-\tin$ and \textit{burn time} as~$\tb\definedas\tf-\tig$. During the coast phase, the vehicle's states passively evolve according to its engine-off dynamics, whereas during the burn phase, the vehicle actively maneuvers in order to achieve its landing objective. Without loss of generality, we define the ignition time epoch as the time at which~$t=0$, where it follows that~$\tc = -\tin$ and that~$\tb = \tf$. This timeline is illustrated in Figure~\ref{fig:traj_timeline}.

The subscripts~$\inertial$ and~$\body$ are used to denote problem parameters expressed in the inertial and body-fixed reference frames~$\iframe$ and~$\bframe$, respectively. We define~$\iframe$ as a surface-fixed Up-East-North coordinate frame whose origin coincides with the landing site. Likewise, we define~$\bframe$ such that its origin coincides with the vehicle's center-of-mass, its $x$-axis points along the vertical axis of the vehicle, its $y$-axis points out the side of the vehicle, and its $z$-axis completes the right-handed system. We use~$\m(t)\in\real_{++}$, $\pos(t)\in\real^3$, and~$\vel(t)\in\real^3$ to respectively denote the mass, position, and velocity states, and~$\thrust(t)\in\real^3$ and~$\aero(t)\in\real^3$ to respectively denote the thrust vector and aerodynamic force. We denote the unit quaternion that parameterizes the transformation from~$\iframe$ to~$\bframe$ by~$\qIB(t)\in\sthree\subset\real^{4}$, and its corresponding direction cosine matrix by~$\cIB(t)\definedas\cIB\big(\qIB(t)\big)\in\sothree$~\cite{Hanson2006}. The conjugate of~$\qIB(t)$ is denoted by~$\qIB^{*}(t)$, where it follows that~$\cBI(t)\definedas\cIB^T(t) = \cIB\big(\qIB^*(t)\big)$. Quaternion multiplication and the identity quaternion are denoted by~$\otimes$ and~$\qidentity$, respectively. The angular velocity of~$\bframe$ relative to~$\iframe$ is denoted by~$\omegaB(t)\in\real^3$, and is expressed in body-fixed coordinates. Lastly, we use~$\bvec{e}_i$ to denote the~$\ith{i}{th}$ basis vector of~$\real^n$, $\crossprod{}{}$ to denote the vector cross product, and~$\;\dotprod{}{}\;$ to denote the vector dot product.

%%%%%%%%%%%%%%%%%%%%%%%%%%%%%%%%%%%%%%%%%%%%%%%%%%%%%%%%%%%%%%%%%%%%%%%%%%%%%%%%%%%%%%

\subsection{Baseline Problem} \label{sec2:6dof}

\subsubsection{Dynamics} As stated in~\sref{sec2:assumptions}, the mass-depletion dynamics are assumed to be an affine function of the thrust magnitude, and are given by
\begin{equation} \label{eq:mass_dyn}
	\mdot(t) = -\mdotalpha\twonorm{\TB(t)}-\mdotbeta\,,
\end{equation}
where~$\mdotalpha\definedas1/\Isp\gstd$ and~$\mdotbeta\definedas\mdotalpha\Pamb\,\Anoz$,~$\Isp$ is the vacuum specific impulse of the engine,~$\gstd$ is standard Earth gravity,~$\Anoz$ is the nozzle exit area of the engine, and~$\Pamb$ is ambient atmospheric pressure. The second term in~\eqref{eq:mass_dyn} represents the specific impulse reduction incurred by atmospheric back-pressure, and assumes that the nozzle remains choked over the allowable throttle range. 

The translational states are governed by the following dynamics
\begin{subequations} \label{eq:trans_dynamics}
	\begin{align}
		\rIdot(t) &= \vI(t)\,, \label{eq:trans_dynamics_a} \\
        \vIdot(t) &= \frac{1}{\m(t)}\FI(t)+\gI, \label{eq:trans_dynamics_b}
	\end{align}
\end{subequations}
where~$\gI\in\real^3$ and~$\FI(t)\definedas\cBI(t)\TB(t) + \aeroI(t)\in\real^{3}$ respectively denote the constant gravitational acceleration and net propulsive and aerodynamic force acting on the vehicle. The thrust vector is expressed in~$\bframe$ coordinates to simplify the attitude dynamics and control constraints that follow. The aerodynamic term~$\aeroI(t)$ is defined~\sref{sec2:aero}.

The attitude states are governed by the following rigid-body attitude dynamics
\begin{subequations} \label{eq:att_dyn}
	\begin{align}
    	\qIBdot(t) &= \frac{1}{2}\skewsymbig\big(\omegaB(t)\big)\qIB(t)\,, \label{eq:att_dyn_a} \\
		\inertia\omegaBdot(t) &= \MB(t) - \crossprod{\omegaB(t)}{\inertia\omegaB(t)}\,, \label{eq:att_dyn_b}
    \end{align}
\end{subequations}
where~$\skewsymbig(\cdot)$ is a skew-symmetric matrix defined such that the quaternion kinematics in~\eqref{eq:att_dyn_a} hold~\cite{Hanson2006}, $\inertia\in\spd{3}$ denotes the body-fixed constant inertia tensor of the vehicle about its center of mass, and~$\MB(t)\definedas\displaystyle{\crossprod{\rTB}{\TB(t)} + \crossprod{\rCPB}{\aeroB(t)}}\in\real^{3}$ denotes the net propulsive and aerodynamic torque acting on the vehicle. The vectors~$\rTB\in\real^3$ and~$\rCPB\in\real^3$ give the constant positions of the engine gimbal pivot point and the aerodynamic center of pressure, respectively. The aerodynamic term $\aeroB(t)$ is defined in~\sref{sec2:aero}.

%%%%%%%%%%%%%%%%%%%%%%%%%%%%%%%%%%%%%%%%%%%%%%%%%%%%%%%%%%%%%%%%%%%%%%%%%%%%%%%%%

\subsubsection{State Constraints} We impose four state constraints in our baseline formulation. First, we constrain the mass of the vehicle to values greater than a minimum dry mass~$\mdry\in\real_{++}$ by enforcing
\begin{equation} \label{eq:minmass}
	\m(t) \geq \mdry.
\end{equation}
Second, we constrain the inertial position to lie inside of a glide-slope cone with half-angle~$\glideslope\in\intei{0\dg}{90\dg}$ and vertex at the origin of~$\iframe$ by enforcing
\begin{equation} \label{eq:glide_slope}
	\dotprod{\ex}{\rI(t)} \geq \tan \glideslope \twonorm{\Hgs \rI(t)},
\end{equation}
where~$\Hgs\definedas[\ey\;\,\ez]^T\in\real^{2 \times 3}$. Third, we define the tilt angle of the vehicle as the angle between the $x$-axes of~$\iframe$ and~$\bframe$, and constrain it to be less than a maximum tilt angle~$\tiltmax\in\intie{0\dg}{90\dg}$ by enforcing
\begin{equation} \label{eq:max_tilt}
	\cos\tiltmax \leq 1-2\twonorm{\Htilt \qIB(t)},
\end{equation}
where~$\Htilt\definedas[\bvec{e}_3\;\,\bvec{e}_4]^T\in\real^{2 \times 4}$ if a scalar-first quaternion convention is used. Fourth, we limit the angular velocity to a maximum value of~$\omegamax\in\real_{++}$ by enforcing
\begin{equation} \label{eq:wmax}
	\twonorm{\omegaB(t)} \leq \omegamax.
\end{equation}

%%%%%%%%%%%%%%%%%%%%%%%%%%%%%%%%%%%%%%%%%%%%%%%%%%%%%%%%%%%%%%%%%%%%%%%%%%%%%%%%%

\subsubsection{Control Constraints} We impose two control constraints in our baseline formulation. First, we impose lower and upper bounds on the thrust magnitude such that
\begin{equation} \label{eq:thrust_bounds}
	0 < \Tmin \leq \twonorm{\TB(t)} \leq \Tmax\,,
\end{equation}
where~$\Tmin$ and~$\Tmax$ are the minimum and maximum allowable thrust magnitudes, respectively. Second, we constrain the thrust vector to lie within a prescribed maximum gimbal angle~$\gimbalmax\in\intii{0\dg}{90\dg}$ relative to the $x$-axis of~$\bframe$ by enforcing
\begin{equation} \label{eq:max_gimbal}
	\cos \gimbalmax \twonorm{\TB(t)} \leq \dotprod{\ex}{\TB(t)}.
\end{equation}

%%%%%%%%%%%%%%%%%%%%%%%%%%%%%%%%%%%%%%%%%%%%%%%%%%%%%%%%%%%%%%%%%%%%%%%%%%%%%%%%%

\subsubsection{Boundary Conditions} We now present a notional set of boundary conditions, with the understanding that they may be modified to accommodate different scenarios. The conditions at the ignition time epoch are given by
\begin{equation} \label{eq:bcs_initial}
	m(\tig) = \mig\,, \quad \rI(\tig) = \rIig\,, \quad \vI(\tig) = \vIig\,, \quad \omegaB(\tig) = \omegaBig\,,
\end{equation}
where $\mig\in\real_{++}$, $\rIig\in\real^3$, and $\vIig\in\real^3$ are the prescribed mass, position, and velocity at ignition time~$\tig$, respectively. We assume that~$\mig>\mdry$. The conditions at the final time epoch are given by
\begin{equation} \label{eq:bcs_final}
	\rI(\tf) = \rIf\,, \quad \vI(\tf) = \vIf\,, \quad \qIB(\tf) = \qIBf\,, \quad \omegaB(\tf) = \omegaBf\,, 
\end{equation}
where~$\vdes\in\real_{+}$ is the prescribed terminal vertical descent speed.

%%%%%%%%%%%%%%%%%%%%%%%%%%%%%%%%%%%%%%%%%%%%%%%%%%%%%%%%%%%%%%%%%%%%%%%%%%%%%%%%%%

\subsection{Free Ignition Time} \label{sec2:fit}

In the baseline problem formulation presented in~\sref{sec2:6dof}, $\rI(\tig)$ and~$\vI(\tig)$ are restricted to prescribed \textit{points} at the time of ignition. We now introduce \textit{the free-ignition-time modification} to~\eqref{eq:bcs_initial}, such that~$\rI(\tig)$ and~$\vI(\tig)$ are constrained to a prescribed \textit{curve} as follows
\begin{equation}
\m(\tig) = \mig\,, \quad \rI(\tig) = \prig(\tc)\,, \quad \vI(\tig) = \pvig(\tc)\,, \quad \omegaB(\tig) = \omegaBig\,,
\label{eq:bcs_initial_fit}
\end{equation}
where the coast time~$\tc\in\intee{0}{\tcmax}$ is included as a solution variable, and~$\prig: \real\rightarrow\real^3$ and~$\pvig: \real\rightarrow\real^3$ are predetermined vector valued polynomials describing an engine-off trajectory. We choose these polynomials to represent an aerodynamics-free trajectory using
\begin{equation}\label{eq:bcs_fit_polynomials}
\prig(\xi) \definedas \rIin + \vIin\,\xi + \frac{1}{2}\gI\xi^2, \quad \pvig(\xi) \definedas \vIin + \gI\xi,
\end{equation}
where~$\rIin$ and~$\vIin$ are prescribed position and velocity vectors at the initial time epoch (see Figure~\ref{fig:traj_timeline}).
Higher order effects (e.g. aerodynamics) can be embedded in~$\prig(\cdot)$ and~$\pvig(\cdot)$ by using higher fidelity models to propagate the vehicle state over a prediction horizon of length~$\tcmax$, and fitting polynomials to the resulting path.

%%%%%%%%%%%%%%%%%%%%%%%%%%%%%%%%%%%%%%%%%%%%%%%%%%%%%%%%%%%%%%%%%%%%%%%%%%%%%%%%%%%

\subsection{Aerodynamic Model} \label{sec2:aero}

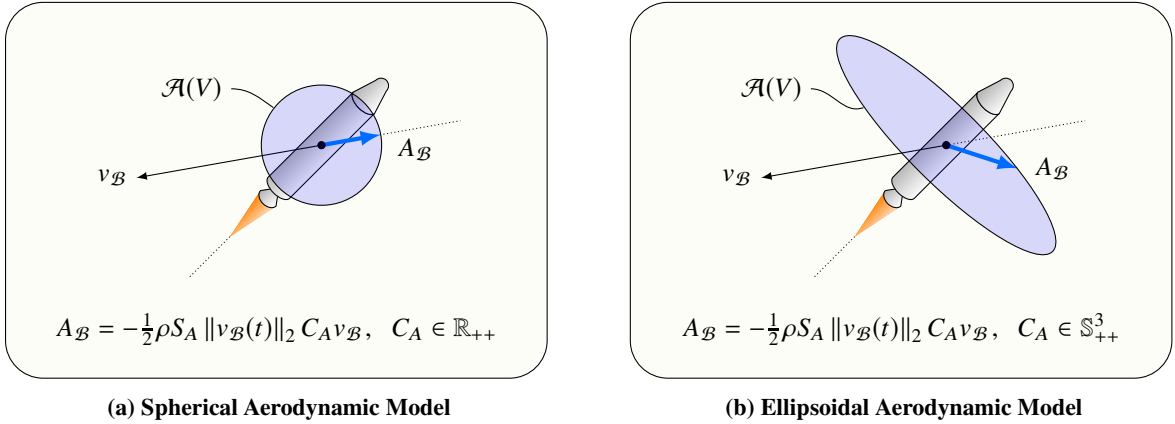
\begin{figure}[t!]
	\begin{subfigure}{0.5\textwidth}
    	\centering
		\begin{tikzpicture}
    \tikzmath{\aa=-45;
    		  \cx=0.6; \cy=0.6;
              \Lcone=1.75; \Acone=30; \Dcone=0.2;
              \labx=-1.2; \laby=0.5;
              \LL=2.5; \LLL=0.80;
              \aaa=-35; \aaaa=\aaa; \fact=0.75; \da=10;
              \lx=sin(\aa)*\LL; \ly=-cos(\aa)*\LL;
              \llx=sin(\aa+\aaa)*\LL; \lly=-cos(\aa+\aaa)*\LL;
              \lllx=-sin(\aa+\aaaa)*\LLL; \llly=cos(\aa+\aaaa)*\LLL;
    }
    
    % Pane %
    \pane{0}{0}{0}{7.2}{5}{0.5}{color=black,fill=beige!20}
    
    % Dotted Lines in Background %
	\draw[black,densely dotted] (\cx,\cy) -- +(\lx,\ly);
    
    % Rocket
    \rocket{\cx}{\cy}{\aa}{1.6}{0.6}{0}
    
    % Velocity Vector %
    \draw[black!100,->] (\cx,\cy) -- +(\llx,\lly);
    \draw[black,densely dotted] (\cx,\cy) -- +(-0.75*\llx,-0.75*\lly);
    \draw (\cx+\llx,\cy+\lly) node[anchor=east] {$\vB$};
    
    % Aero Vector %
    \draw[acol,ultra thick,->] (\cx,\cy) -- +(\lllx,\llly);
    \draw[black!100,fill=black] (\cx,\cy) circle (0.3ex);
    \draw (\cx+\lllx,\cy+\llly-0.2) node[anchor=west] {$\;\aeroB$};
    
    % Ellipse %
    \fill[color=blue!100,opacity=0.15] (\cx,\cy) ellipse (0.80 and 0.80);
%    \shade[rotate around={\aa:(\cx,\cy)},ball color=blue,fill opacity=0.2] (\cx,\cy) ellipse (0.80 and 0.80);
    \draw[rotate around={\aa:(\cx,\cy)}] (\cx,\cy) ellipse (0.80 and 0.80);
    \draw [black!100] (\cx-0.64,\cy+0.5) to[out=145,in=0] (\cx-1.2,\cy+0.75) node[anchor=east,color=black!100] {$\setaero$};
    
%    \draw(0,-1.4) node[anchor=north] {$\aeroB = -\displaystyle{\frac{1}{2}}\density\Sa\twonorm{\vB(t)}\Ca\vB\,,\;\;\Ca\in\real_{++} $};
    \draw(0,-1.5) node[anchor=north] {$\aeroB = -\frac{1}{2}\density\Sa\twonorm{\vB(t)}\Ca\vB\,,\;\;\Ca\in\real_{++} $};
\end{tikzpicture}
    	\caption{Spherical Aerodynamic Model}
    	\label{fig:aerodynamic_model_a}
   	\end{subfigure}
    \begin{subfigure}{0.5\textwidth}
    	\centering
		\begin{tikzpicture}
    \tikzmath{\aa=-45;
    		  \cx=0.6; \cy=0.6;
              \Lcone=1.75; \Acone=30; \Dcone=0.2;
              \labx=-1.2; \laby=0.5;
              \LL=2.5; \LLL=1.0;
              \aaa=-35; \aaaa=-63; \fact=0.75; \da=10;
              \lx=sin(\aa)*\LL; \ly=-cos(\aa)*\LL;
              \llx=sin(\aa+\aaa)*\LL; \lly=-cos(\aa+\aaa)*\LL;
              \lllx=-sin(\aa+\aaaa)*\LLL; \llly=cos(\aa+\aaaa)*\LLL;
    }
    
    % Pane %
    \pane{0}{0}{0}{7.2}{5}{0.5}{color=black,fill=beige!20}
    
    % Dotted Lines in Background %
	\draw[black,densely dotted] (\cx,\cy) -- +(\lx,\ly);
    
    % Rocket
    \rocket{\cx}{\cy}{\aa}{1.6}{0.6}{0}
    
    % Velocity Vector %
    \draw[black!100,->] (\cx,\cy) -- +(\llx,\lly);
    \draw[black,densely dotted] (\cx,\cy) -- +(-0.75*\llx,-0.75*\lly);
    \draw (\cx+\llx,\cy+\lly) node[anchor=east] {$\vB$};
    
    % Aero Vector %
    \draw[acol,ultra thick,->] (\cx,\cy) -- +(\lllx,\llly);
    \draw[black!100,fill=black] (\cx,\cy) circle (0.3ex);
    \draw (\cx+\lllx,\cy+\llly) node[anchor=west] {$\;\aeroB$};
    
    % Ellipse %
    \fill[rotate around={\aa:(\cx,\cy)},color=blue!100,fill opacity=0.15] (\cx,\cy) ellipse (2 and 0.5);
%    \shade[rotate around={\aa:(\cx,\cy)},ball color=blue,fill opacity=0.2] (\cx,\cy) ellipse (2 and 0.5);
    \draw[rotate around={\aa:(\cx,\cy)}] (\cx,\cy) ellipse (2 and 0.5);
    \draw [black!100] (\cx-1.1,\cy+0.5) to[out=200,in=0] (\cx-1.8,\cy+0.75) node[anchor=east,color=black!100] {$\setaero$};
    
%    \draw(0,-1.4) node[anchor=north] {$\aeroB = -\displaystyle{\frac{1}{2}}\density\Sa\twonorm{\vB(t)}\Ca\vB\,,\;\;\Ca\in\spd{3} $};
    \draw(0,-1.5) node[anchor=north] {$\aeroB = -\frac{1}{2}\density\Sa\twonorm{\vB(t)}\Ca\vB\,,\;\;\Ca\in\spd{3} $};
\end{tikzpicture}
    	\caption{Ellipsoidal Aerodynamic Model}
    	\label{fig:aerodynamic_model_b}
   	\end{subfigure}
    \caption{Depictions of the aerodynamic models introduced in~\sref{sec2:aero}. The spherical model generates drag forces only, whereas the ellipsoidal model can generate both lift and drag forces. The ellipsoidal model is a contribution of this paper, and is a generalization of the spherical model, which was introduced in~\cite{szmuk2016successive}.}
    \label{fig:aerodynamic_model}
\end{figure}

We now introduce a tractable aerodynamic model that approximates the relationship between the aerodynamic force and the velocity vector. The model expresses the aerodynamic force in~$\bframe$ coordinates as follows
\begin{equation} \label{eq:dragB}
	\aeroB(t) = - \frac{1}{2}\density V(t)\Sa\Ca\cBI(t)\vI(t)\,,\quad\Ca\in\spd{3}\,,
\end{equation}
where~$\density$ is the ambient atmospheric density, $V(t)\definedas\twonorm{\vI(t)}$, and~$\Sa\in\real_{++}$ is a constant reference area. We refer to~$\Ca$ as the aerodynamic coefficient matrix, and emphasize that it is a symmetric-positive-definite matrix that does not conform to the standard scalar definition. Our definition of~$\Ca$ ensures that for any~$V\in\real_{+}$ the set
$\setaero\definedas\big\{ \aeroB : \aeroB = -\frac{1}{2}\density V^2\Sa\Ca\hat{\bvec{v}}_\body\,,\; \hat{\bvec{v}}_\body\in\mathcal{S}^2 \,\big\}$ defines a fixed ellipsoid in~$\bframe$ coordinates. Since most rocket-powered vehicles are approximately axisymmetric, we assume~$\Ca = \diag{[c_{a,x},\; c_{a,yz},\; c_{a,yz}]}$, where~$c_{a,x}$ and~$c_{a,yz}$ are positive scalars. This assumption aligns the principal axes of~$\setaero$ with the axes of~$\bframe$.

If~$c_{a,x} = c_{a,yz}$, then~$\aeroB(t)$ is always anti-parallel to~$\vB(t)$. In this case,~$\aeroB(t)$ may be interpreted as a pure drag force, and the model recovers the aerodynamic drag model used in~\cite{szmuk2016successive}. Since the set~$\setaero$ corresponding to this choice of~$\Ca$ defines a sphere, we refer to the corresponding model as the \textit{spherical aerodynamic model}, illustrated in Figure~\ref{fig:aerodynamic_model_a}. Under the assumptions of the spherical model, the product~$\cBI(t) \Ca \cIB(t)$ simplifies to~$c_{a,x} I_{3 \times 3}$, thus rendering~$\aeroI(t) = \cBI(t)\aeroB(t)$ independent of attitude.

Alternatively, if~$c_{a_x} \neq c_{a,yz}$, then~$\aeroB(t)$ can also have components orthogonal to~$\vB(t)$. In this case,~$\aeroB(t)$ may be interpreted as the vector sum of a drag and lift force. Furthermore, if we assume that~$c_{a,x} < c_{a,yz}$, we ensure that the vehicle experiences minimum drag when~$\vB$ is aligned with the $x$-axis of~$\bframe$, and that the lift component of $\aeroB$ points in the correct direction. Since the set~$\setaero$ corresponding to this choice of~$\Ca$ defines an ellipsoid, we refer to the corresponding model as the \textit{ellipsoidal aerodynamic model}, illustrated in Figure~\ref{fig:aerodynamic_model_b}.

% /* LEAVE THIS COMMENT IN 
%The force $\aeroB(t)$ can be equivalently expressed in $\iframe$ coordinates as
%
%\begin{equation} \label{eq:dragI}
%	\aeroI(t) = \cBI(t)\aeroB(t) = -\frac{1}{2}\density\Sa\twonorm{\vI(t)}\CaI(t)\vI(t)\,,\quad\Ca\in\spd{3}\,,
%\end{equation}
%
%where
%
%\begin{equation} \label{eq:coefdragI}
%	\CaI(t)\definedas\cIB^{T}(t)\Ca\cIB(t).
%\end{equation}
%
%Note that when $\Ca$ is a multiple of the identity matrix (i.e. for the spherical model), then $\CaI(t) = \Ca$.
% LEAVE THIS COMMENT IN */

%%%%%%%%%%%%%%%%%%%%%%%%%%%%%%%%%%%%%%%%%%%%%%%%%%%%%%%%%%%%%%%%%%%%%%%%%%%

\subsection{State-Triggered Constraints} \label{sec2:stc}

In this section, we introduce the most important contribution of this paper: a continuous formulation of \textit{state-triggered constraints} (STCs). The most common type of constraints seen in the optimal control literature are enforced over predetermined time intervals; we refer to such constraints as \textit{temporally-scheduled constraints}. In contrast, an STC is enforced only when a state-dependent condition is satisfied, and emulates a constraint gated by an \textit{if}-statement conditioned on the solution variables. Thus, an optimal control problem containing an STC determines its solution variables with a simultaneous understanding of how the constraint affects the optimization, and of how the optimization enables or disables the constraint. While the resulting continuous formulation is still non-convex, we have found that it is readily handled by the successive convexification framework \cite{SCvx_2016arXiv,SCvx_2017arXiv,SCvx_cdc16}, as demonstrated in~\sref{sec:numerical_results}.

%%%%%%%%%%%%%%%%%%%%%%%%%%%%%%%%%%%%%%%%%%%%%%%%%%%%%%%%%%%%%%%%%%%%%%%%%%%%%%%%%

\subsubsection{Formal Definition of State-Triggered Constraints} \label{sec2:stc_formal_def}

Formally, we define a \textit{state-triggered constraint} as an equality constraint that is enforced conditionally according to the following logical statement
\begin{equation} \label{eq:stc_orig}
	\stcg(\stcx) < 0 \;\Rightarrow\; \stcf(\stcx) = 0\,,
\end{equation}
where~$\stcx\in\real^\stcnx$ represents the optimization variable of the parent problem,~$\stcg(\cdot) : \real^\stcnx \rightarrow \real$ is a piecewise continuously differentiable function called the \textit{trigger function}, and~$\stcf(\cdot) : \real^\stcnx \rightarrow \real$ is a piecewise continuously differentiable function called the \textit{constraint function}. Accordingly, we refer to~$\stcg(\stcx)<0$ as the \textit{trigger condition}, and~$\stcf(\stcx) = 0$ as the \textit{constraint condition}.

The logical implication in~\eqref{eq:stc_orig} states that if the trigger condition is satisfied, then the constraint condition is enforced. By De Morgan's Law,~\eqref{eq:stc_orig} also implies satisfaction of the contrapositive. However, we emphasize that the satisfaction of the constraint condition does not imply satisfaction of the trigger condition (see~\sref{sec:numerical_results}). 

\begin{remark} \label{rem:stc_eq_ineq}
	The constraint condition in~\eqref{eq:stc_orig} can be used to represent an inequality constraint by augmenting~$\stcx$ with a non-negative slack variable and modifying~$\stcf(\stcx)$ accordingly~\cite{BoydConvex}.
\end{remark}

%%%%%%%%%%%%%%%%%%%%%%%%%%%%%%%%%%%%%%%%%%%%%%%%%%%%%%%%%%%%%%%%%%%%%%%%%%%%%%%%%

\subsubsection{Continuous State-Triggered Constraints} \label{sec2:stc_cont_form}

The STC expressed in~\eqref{eq:stc_orig} represents a binary decision that is not readily incorporated into a continuous optimization framework. We address this issue by introducing \textit{continuous state-triggered constraints} (cSTCs), which represent the logical implication in~\eqref{eq:stc_orig} using the auxiliary variable~$\stcs\in\real_{+}$ and the system of equations
\begin{subequations} \label{eq:stc_cont}
	\begin{align}
    	\stcg(\stcx) + \stcs &\geq 0\,, \label{eq:stc_cont_a} \\
        \stcs &\geq 0\,, \label{eq:stc_cont_b} \\
        \stcs\cdot\stcf(\stcx) &= 0. \label{eq:stc_cont_c}
    \end{align}
\end{subequations}
The geometry of the cSTC in~\eqref{eq:stc_cont} is shown on the left side of Figure~\ref{fig:cstc_geom}, where the lower-left axes show the feasible set of the STC in~\eqref{eq:stc_orig}. The formulation in~\eqref{eq:stc_cont} ensures that~$\stcs$ is strictly positive if the trigger condition is satisfied. Since~$\stcs>0$ implies that~\eqref{eq:stc_cont_c} holds if and only if~$\stcf(\stcx) = 0$, \eqref{eq:stc_orig} and~\eqref{eq:stc_cont} are logically equivalent. Subject to mild assumptions on~$\stcg(\cdot)$ and~$\stcf(\cdot)$,~\eqref{eq:stc_cont} admits a solution for any~$\stcx\in\real^{\stcnx}$. 

%%%%%%%%%%%%%%%%%%%%%%%%%%%%%%%%%%%%%%%%%%%%%%%%%%%%%%%%%%%%%%%%%%%%%%%%%%%%%%%%%

\subsubsection{Improved Formulation Using Linear Complementarity} \label{sec2:stc_lcp}

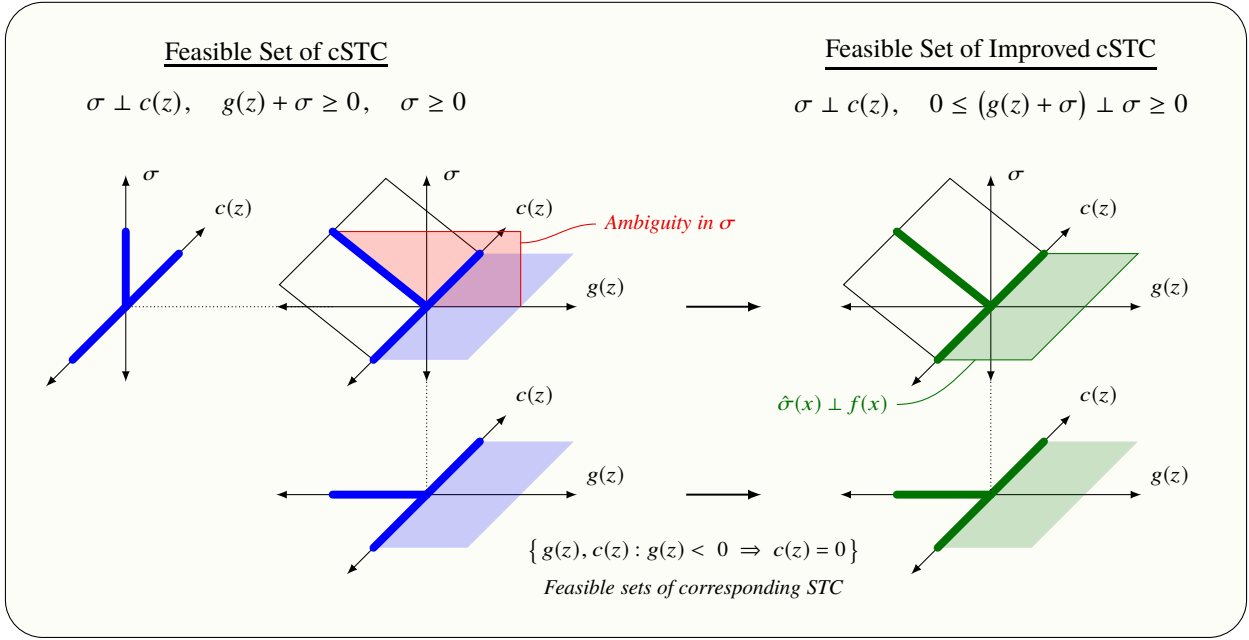
\begin{figure}[t!]
	\centering
  	\begin{tikzpicture}
	\tikzmath{\xa=0; \ya=45; \Lx=2.00; \Ly=1.5; \Lzt=1.75; \Lzb=1.00; \L2=2.50; \L3=4.0; \L4=7.5;
			  \cx1=0; \cy1=0;
              \cx2=\cx1; \cy2=\cy1-\L2;
              \cx3=\cx1-cos(\xa)*\L3; \cy3=\cy1+sin(\xa)*\L3;
              \cx4=\cx1+cos(\xa)*\L4; \cy4=\cy1-sin(\xa)*\L4;
              \cx5=\cx4; \cy5=\cy4-\L2;
							\cx6=0.5*\cx1+0.5*\cx3; \cy6=0.5*\cy1+0.5*\cy3;
              \LL=1; \LLX=1.25*\LL; \LLYx=sin(\ya)*\LL; \LLYy=cos(\ya)*\LL;
              \txHt1=1.60; \txHt2=1.25;
							\darkness=30;
    }
    
    % Panes %
		\cpane{\cx3-1.6}{\cx4+3.4}{\cy3+\Lzt+\txHt1+0.70}{\cy2-\Lzb-0.90}{0}{0.5}{color=black,fill=beige!20}
		
		% Titles %
    \draw (\cx6,\cy6+\Lzt+\txHt1) node[anchor=center] {\underline{Feasible Set of cSTC}};
    \draw (\cx4,\cy4+\Lzt+\txHt1) node[anchor=center] {\underline{Feasible Set of Improved cSTC}};
    
		% Left subtitle %
		\draw (\cx6,\cy6+\Lzt+\txHt2) node[anchor=north] {$\stcs \perp \stcf(\stcx)\,,\quad \stcg(\stcx)+\stcs \geq 0\,,\quad \stcs \geq 0$};
		
		%\munderbrace{\cx6-1.95}{\cy6+\Lzt+\txHt2-0.6}{0}{1.40}{0.0}{0.1}{black!\darkness}
		%\draw [black!\darkness] (\cx6-1.95,\cy6+\Lzt+\txHt2-0.6) to[out=-90,in=70] (\cx3-0.6*\LLYx,\cy3+1.2) to[out=-110,in=180] (\cx3,\cy3+0.5*\LL);%(\cx3-0.6*\LLYx,\cy3-0.6*\LLYy);
		
		%\munderbrace{\cx6-0.30}{\cy6+\Lzt+\txHt2-0.6}{0}{1.00}{0.0}{0.1}{black!\darkness}
		%\draw [black!\darkness] (\cx6-0.30,\cy6+\Lzt+\txHt2-0.6)
		%	to[out=-90,in=180-\ya] (\cx1-\LLX-0.2*\LLYx,\cy1+\LL-0.2*\LLYy)
		%	-- (\cx1-0.2*\LLYx,\cy1-0.2*\LLYy);
		
		%\munderbrace{\cx6+1.61}{\cy6+\Lzt+\txHt2-0.6}{0}{2.00}{0.0}{0.1}{black!\darkness}
		%\draw [black!\darkness] (\cx6+1.61,\cy6+\Lzt+\txHt2-0.6) to[out=-90,in=90] (\cx1-\LLX+\LLYx,\cy1+\LL+\LLYy);
		
		% Right subtitle %
    \draw (\cx4,\cy4+\Lzt+\txHt2) node[anchor=north] {$\stcs \perp \stcf(\stcx)\,,\quad 0 \leq \big(\stcg(\stcx)+\stcs\big) \perp \stcs \geq 0$};
		
    % Axes %
    \threeaxes      {\cx1}{\cy1}{\Lx}{\Ly}{\Lzt}{\Lzb}{\xa}{\ya}
    \threeaxeslabelx{\cx1}{\cy1}{\Lx}{\Ly}{\Lzt}{\Lzb}{\xa}{\ya}{\footnotesize $\stcg(\stcx)$}
    \threeaxeslabely{\cx1}{\cy1}{\Lx}{\Ly}{\Lzt}{\Lzb}{\xa}{\ya}{\footnotesize $\stcf(\stcx)$}
    \threeaxeslabelz{\cx1}{\cy1}{\Lx}{\Ly}{\Lzt}{\Lzb}{\xa}{\ya}{\footnotesize $\stcs$}
    
    \threeaxes      {\cx2}{\cy2}{\Lx}{\Ly}{0   }{0   }{\xa}{\ya}
    \threeaxeslabelx{\cx2}{\cy2}{\Lx}{\Ly}{0   }{0   }{\xa}{\ya}{\footnotesize $\stcg(\stcx)$}
    \threeaxeslabely{\cx2}{\cy2}{\Lx}{\Ly}{0   }{0   }{\xa}{\ya}{\footnotesize $\stcf(\stcx)$}
    
    \threeaxes      {\cx3}{\cy3}{0  }{\Ly}{\Lzt}{\Lzb}{\xa}{\ya}
    \threeaxeslabely{\cx3}{\cy3}{0  }{\Ly}{\Lzt}{\Lzb}{\xa}{\ya}{\footnotesize $\stcf(\stcx)$}
    \threeaxeslabelz{\cx3}{\cy3}{0  }{\Ly}{\Lzt}{\Lzb}{\xa}{\ya}{\footnotesize $\stcs$}
    
    \threeaxes{\cx4}{\cy4}{\Lx}{\Ly}{\Lzt}{\Lzb}{\xa}{\ya}
    \threeaxeslabelx{\cx4}{\cy4}{\Lx}{\Ly}{\Lzt}{\Lzb}{\xa}{\ya}{\footnotesize $\stcg(\stcx)$}
    \threeaxeslabely{\cx4}{\cy4}{\Lx}{\Ly}{\Lzt}{\Lzb}{\xa}{\ya}{\footnotesize $\stcf(\stcx)$}
    \threeaxeslabelz{\cx4}{\cy4}{\Lx}{\Ly}{\Lzt}{\Lzb}{\xa}{\ya}{\footnotesize $\stcs$}
    
    \threeaxes      {\cx5}{\cy5}{\Lx}{\Ly}{0   }{0   }{\xa}{\ya}
    \threeaxeslabelx{\cx5}{\cy5}{\Lx}{\Ly}{0   }{0   }{\xa}{\ya}{\footnotesize $\stcg(\stcx)$}
    \threeaxeslabely{\cx5}{\cy5}{\Lx}{\Ly}{0   }{0   }{\xa}{\ya}{\footnotesize $\stcf(\stcx)$}
    
    % cSTC %
    \draw[color=black!100] (\cx1+\LLYx,\cy1+\LLYy)
                        -- (\cx1+\LLYx-\LLX,\cy1+\LLYy+\LL)
                        -- (\cx1-\LLYx-\LLX,\cy1-\LLYy+\LL)
                        -- (\cx1-\LLYx,\cy1-\LLYy)
                        -- cycle;
%    \draw[color=blue!100]  (\cx1+\LLYx,\cy1+\LLYy)
%                        -- (\cx1+\LLYx+\LLX,\cy1+\LLYy)
%                        -- (\cx1-\LLYx+\LLX,\cy1-\LLYy)
%                        -- (\cx1-\LLYx,\cy1-\LLYy)
%                        -- cycle;
    \fill[color=blue!100,opacity=0.2] (\cx1+\LLYx,\cy1+\LLYy)
                                   -- (\cx1+\LLYx+\LLX,\cy1+\LLYy)
                                   -- (\cx1-\LLYx+\LLX,\cy1-\LLYy)
                                   -- (\cx1-\LLYx,\cy1-\LLYy)
                                   -- cycle;
    \fill[color=red,opacity=0.2]  (\cx1-\LLX,\cy1+\LL)
    						   -- (\cx1+\LLX,\cy1+\LL)
                               -- (\cx1+\LLX,\cy1)
                               -- (\cx1,\cy1)
                               -- cycle;
    \draw[color=darkred]  (\cx1-\LLX,\cy1+\LL)
    				   -- (\cx1+\LLX,\cy1+\LL)
                       -- (\cx1+\LLX,\cy1);
    \draw[blue!100,line width=3pt,line cap=round] (\cx1,\cy1) -- +(-\LLX,\LL);
    \draw[blue!100,line width=3pt,line cap=round] (\cx1,\cy1) -- +( \LLYx, \LLYy);
    \draw[blue!100,line width=3pt,line cap=round] (\cx1,\cy1) -- +(-\LLYx,-\LLYy);
    
    % STC Feasible set 1 %
    \fill[color=blue!100,opacity=0.2] (\cx2+\LLYx,\cy2+\LLYy)
                                   -- (\cx2+\LLYx+\LLX,\cy2+\LLYy)
                                   -- (\cx2-\LLYx+\LLX,\cy2-\LLYy)
                                   -- (\cx2-\LLYx,\cy2-\LLYy)
                                   -- cycle;
%    \draw[color=blue!100]  (\cx2+\LLYx,\cy2+\LLYy)
%                        -- (\cx2+\LLYx+\LLX,\cy2+\LLYy)
%                        -- (\cx2-\LLYx+\LLX,\cy2-\LLYy)
%                        -- (\cx2-\LLYx,\cy2-\LLYy)
%                        -- cycle;
    \draw[blue!100,line width=3pt,line cap=round] (\cx2,\cy2) -- +(-\LLX,0);
    \draw[blue!100,line width=3pt,line cap=round] (\cx2,\cy2) -- +( \LLYx, \LLYy);
    \draw[blue!100,line width=3pt,line cap=round] (\cx2,\cy2) -- +(-\LLYx,-\LLYy);
    
    % Complementarity constraint %
    \draw[blue!100,line width=3pt,line cap=round] (\cx3,\cy3) -- +(0,\LL);
    \draw[blue!100,line width=3pt,line cap=round] (\cx3,\cy3) -- +( \LLYx, \LLYy);
    \draw[blue!100,line width=3pt,line cap=round] (\cx3,\cy3) -- +(-\LLYx,-\LLYy);
    
    % Improved cSTC %
    \draw[color=black!100] (\cx4+\LLYx,\cy4+\LLYy)
                        -- (\cx4+\LLYx-\LLX,\cy4+\LLYy+\LL)
                        -- (\cx4-\LLYx-\LLX,\cy4-\LLYy+\LL)
                        -- (\cx4-\LLYx,\cy4-\LLYy)
                        -- cycle;
    \draw[color=darkgreen!100] (\cx4+\LLYx,\cy4+\LLYy)
                        -- (\cx4+\LLYx+\LLX,\cy4+\LLYy)
                        -- (\cx4-\LLYx+\LLX,\cy4-\LLYy)
                        -- (\cx4-\LLYx,\cy4-\LLYy)
                        -- cycle;
    \fill[color=darkgreen!100,opacity=0.2] (\cx4+\LLYx,\cy4+\LLYy)
                                        -- (\cx4+\LLYx+\LLX,\cy4+\LLYy)
                                        -- (\cx4-\LLYx+\LLX,\cy4-\LLYy)
                                        -- (\cx4-\LLYx,\cy1-\LLYy)
                                        -- cycle;
    \draw[darkgreen!100,line width=3pt,line cap=round] (\cx4,\cy4) -- +(-\LLX,\LL);
    \draw[darkgreen!100,line width=3pt,line cap=round] (\cx4,\cy4) -- +( \LLYx, \LLYy);
    \draw[darkgreen!100,line width=3pt,line cap=round] (\cx4,\cy4) -- +(-\LLYx,-\LLYy);
    
    % STC Feasible set 2 %
    \fill[color=darkgreen!100,opacity=0.2] (\cx5+\LLYx,\cy5+\LLYy)
                                        -- (\cx5+\LLYx+\LLX,\cy5+\LLYy)
                                        -- (\cx5-\LLYx+\LLX,\cy5-\LLYy)
                                        -- (\cx5-\LLYx,\cy5-\LLYy)
                                        -- cycle;
    \draw[darkgreen!100,line width=3pt,line cap=round] (\cx5,\cy5) -- +(-\LLX,0);
    \draw[darkgreen!100,line width=3pt,line cap=round] (\cx5,\cy5) -- +( \LLYx, \LLYy);
    \draw[darkgreen!100,line width=3pt,line cap=round] (\cx5,\cy5) -- +(-\LLYx,-\LLYy);
%    \draw[color=darkgreen!100] (\cx5+\LLYx,\cy5+\LLYy)
%                        -- (\cx5+\LLYx+\LLX,\cy5+\LLYy)
%                        -- (\cx5-\LLYx+\LLX,\cy5-\LLYy)
%                        -- (\cx5-\LLYx,\cy5-\LLYy)
%                        -- cycle;
    
    % Dotted lines %
    \draw[color=black,densely dotted] (\cx1-\LLX,\cy1) -- (\cx3,\cy3);
    \draw[color=black,densely dotted] (\cx1,\cy1-\Lzb) -- (\cx2,\cy2);
    \draw[color=black,densely dotted] (\cx4,\cy4-\Lzb) -- (\cx5,\cy5);
    
    % Arrows %
    \draw [black,line width=0.75pt,->] (0.5*\cx1+0.5*\cx4-0.50+0.2,0.5*\cy1+0.5*\cy4) -- +(1.0,0);
    \draw [black,line width=0.75pt,->] (0.5*\cx1+0.5*\cx4-0.50+0.2,0.5*\cy2+0.5*\cy5) -- +(1.0,0);
    
    % Primary labels %
    \draw [darkred] (\cx1+\LLX,\cy1+0.85*\LL) to[out=0,in=-180] (\cx1+2.25,\cy1+1.1) node[anchor=west,color=darkred] {\footnotesize \textit{Ambiguity in} $\stcs$};
		\draw [darkgreen] (\cx4-\LLYx+0.4*\LLX,\cy4-\LLYy) %		-- (\cx4-\LLX-0.3*\LLYx,\cy4+\LL-0.3*\LLYy)
		to[out=-90-\ya,in=0] (\cx4-1.25,\cy4-1.30) node[anchor=east,color=darkgreen] {\footnotesize $\stcshat(\bvec{x}) \perp f(\bvec{x})$};
    
    % Bottom labels %
    \draw (0.5*\cx2+0.5*\cx5-0.2,\cy2-\Lzb+0.25) node[anchor=center] {\footnotesize $\big\{\,\stcg(\stcx)\,,\stcf(\stcx) : \stcg(\stcx) < \;0 \;\Rightarrow\; \stcf(\stcx)=0\,\big\}$};
    \draw (0.5*\cx2+0.5*\cx5-0.2,\cy2-\Lzb-0.25) node[anchor=center] {\footnotesize \textit{Feasible sets of corresponding STC}};
		%\munderbrace{0.5*\cx2+0.5*\cx5}{\cy2-\Lzb-0.45}{0}{10.5}{0.15}{0.15}{black}
		
		%\draw [darkgreen] (\cx5-\LLYx+0.75*\LLX,\cy5-\LLYy) to[out=-90-\ya,in=0] (\cx5-1.4,\cy5-1.40);
		%\draw [blue] (\cx2-\LLYx+0.75*\LLX,\cy2-\LLYy) -- (\cx2-1.3*\LLYx+0.75*\LLX,\cy2-1.3*\LLYy) to[out=-90-\ya,in=-180] (\cx2+1.1,\cy2-1.40);
		%\draw [blue] (\cx2-\LLYx+0.5*\LLX,\cy2-\LLYy) to[out=-90+\ya,in=-180] (\cx2+1.1,\cy2-1.40);

\end{tikzpicture}
  	\caption{A geometrical interpretation of cSTCs. The red and blue sets on the left represent the geometry of~\eqref{eq:stc_cont}, whereas the green sets on the right represent the geometry of~\eqref{eq:lcp_ad}. The bottom two axes show that the feasible sets of both cSTC formulations comply with that of the STC in~\eqref{eq:stc_orig}, despite the removal of the ambiguity in $\stcs$ in the improved formulation.}
  	\label{fig:cstc_geom}
\end{figure}

As illustrated in Figure~\ref{fig:cstc_geom},~\eqref{eq:stc_cont} does not admit a unique~$\stcs$ given $\stcx$. To resolve this ambiguity, we augment~\eqref{eq:stc_cont_a} and~\eqref{eq:stc_cont_b} to form a complementarity constraint and obtain the following improved cSTC formulation
\begin{subequations} \label{eq:lcp_ad}
  \begin{align}
      0 \leq \stcs &\perp \big(\stcg(\stcx)+\stcs\big) \geq 0\,, \label{eq:lcp_a} \\
      \stcs &\perp \stcf(\stcx)\,,\label{eq:lcp_d}
  \end{align}
\end{subequations}
where~$a \perp b$ is used to denote~$a \cdot b = 0$. For a given~$\stcx$,~\eqref{eq:lcp_a} forms a linear complementarity problem (LCP) in~$\stcs$~\cite{Cottle1992}. This problem has a unique solution~$\stcshat$ that varies continuously in~$\stcg(\stcx)$~\cite{Cottle1992}, and therefore in~$\stcx$. The analytical solution to the LCP is given by
\begin{equation}
	\stcshat(\stcx) \definedas -\min\big(\stcg(\stcx)\,,0\big). \label{eq:lcp_sol}
\end{equation}
Substituting~\eqref{eq:lcp_sol} into~\eqref{eq:lcp_ad} guarantees satisfaction of
~\eqref{eq:lcp_a}, resulting in the following equality constraint
\begin{equation} \label{eq:stc_lcp}
	\stch(\stcx) \definedas -\min\big(\stcg(\stcx)\,,0\big)\cdot\stcf(\stcx) = 0\,,
\end{equation}
where the negative sign is retained to allow~\eqref{eq:stc_lcp} to be formulated as an inequality (see Remark~\ref{rem:stc_eq_ineq}). Thus, the ambiguity in~$\stcs$ is resolved by replacing the constraints in~\eqref{eq:stc_cont} with the logically equivalent constraint given in~\eqref{eq:stc_lcp}, which complies with the feasible set of the corresponding STC defined in~\eqref{eq:stc_orig}. This can be seen in the two rightmost axes of Figure~\ref{fig:cstc_geom}. Subsequent sections use the improved cSTC formulation in lieu of the original one.

%%%%%%%%%%%%%%%%%%%%%%%%%%%%%%%%%%%%%%%%%%%%%%%%%%%%%%%%%%%%%%%%%%%%%%%%%%%%%%%%%

\subsubsection{Example Application} \label{sec2:stc_ex_app}

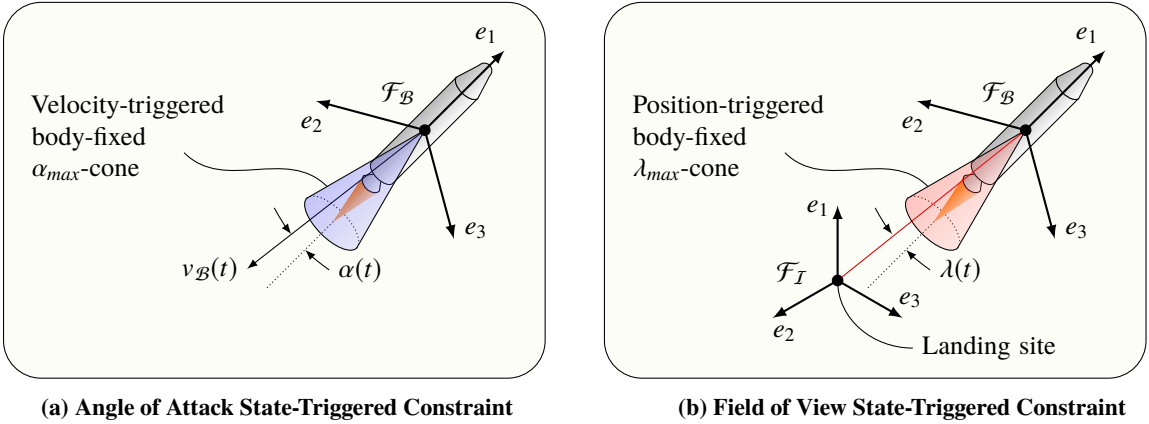
\begin{figure}[tb]
	\begin{subfigure}{0.5\textwidth}
    	\centering
		\begin{tikzpicture}
    \tikzmath{\aa=-45;
    		  \cx=2; \cy=0.8;
              \Lcone=1.75; \Acone=30; \Dcone=0.2;
              \labx=-1.2; \laby=0.5;
              \LL=3.0; \aaa=-7; \fact=0.75; \da=10;
              \lx=sin(\aa)*\LL; \ly=-cos(\aa)*\LL;
              \llx=sin(\aa+\aaa)*\LL; \lly=-cos(\aa+\aaa)*\LL;
              \lllx=sin(\aa-0.5*\Acone)*0.5*\LL; \llly=-cos(\aa-0.5*\Acone)*0.5*\LL;
    }
    
    % Pane %
    \pane{0}{0}{0}{7.2}{5}{0.5}{color=black,fill=beige!20}
    
    % Dotted Lines in Background %
	\draw[black,densely dotted] (\cx,\cy) -- +(\lx,\ly);
    \coneback{\cx}{\cy}{\aa}{\Lcone}{\Acone}{\Dcone}{color=black!100,densely dotted}
    
    % Rocket
    \rocket{\cx}{\cy}{\aa}{1.6}{0.6}{0}
    
    % Velocity Vector %
    \draw[black!100,->] (\cx,\cy) -- +(\llx,\lly);
    \draw (\cx+\llx,\cy+\lly) node[anchor=east] {$\vB(t)$};
    
    % AoA Cone %
    \cone{\cx}{\cy}{\aa}{\Lcone}{\Acone}{\Dcone}{color=blue!100,opacity=0.10,shading=axis,shading angle=90,left color=blue!100,right color=blue!0};
    \draw (\labx+0.85,\laby+0.2) node[text width=2.75cm,align=left,anchor=east] {Velocity-triggered \\ body-fixed \\ $\aoamax$-cone};
    \draw (\cx+\lllx,\cy+\llly) to[out=135,in=-45] (\labx,\laby);
    
    % Axes %
    \isoaxes{\cx}{\cy}{90+\aa}{1.5}{12}{$\bframe$}{$\bvec{e}_1$}{$\bvec{e}_2$}{$\bvec{e}_3$}
    
    % AoA Dimension %
    \draw[color=black!100,<-] (\cx+\fact*\lx,\cy+\fact*\ly) arc(-90+\aa:-90+\aa+\da:2*\fact*\LL) node[anchor=west] {$\aoa(t)$};
    \draw[color=black!100,<-] (\cx+\fact*\llx,\cy+\fact*\lly) arc(-90+\aa+\aaa:-90+\aa+\aaa-\da:2*\fact*\LL);
\end{tikzpicture}
    	\caption{Angle of Attack State-Triggered Constraint}
    	\label{fig:stc_ex_aoa}
   	\end{subfigure}
    \begin{subfigure}{0.5\textwidth}
    	\centering
		\begin{tikzpicture}
    \tikzmath{\aa=-45;
    		  \cx=2; \cy=0.8;
              \Lcone=1.75; \Acone=30; \Dcone=0.2;
              \labx=-1.2; \laby=0.5;
              \LL=3.0; \aaa=-7; \fact=0.75; \da=10;
              \lx=sin(\aa)*\LL; \ly=-cos(\aa)*\LL;
              \llx=sin(\aa+\aaa)*\LL; \lly=-cos(\aa+\aaa)*\LL;
              \lllx=sin(\aa-0.5*\Acone)*0.5*\LL; \llly=-cos(\aa-0.5*\Acone)*0.5*\LL;
              \cix=-0.5; \ciy=-1.2;
    }
    
    % Pane %
    \pane{0}{0}{0}{7.2}{5}{0.5}{color=black,fill=beige!20}
    
    % Dotted Lines in Background %
	\draw[black,densely dotted] (\cx,\cy) -- +(\lx,\ly);
    \coneback{\cx}{\cy}{\aa}{\Lcone}{\Acone}{\Dcone}{color=black!100,densely dotted}
    
    % Rocket
    \rocket{\cx}{\cy}{\aa}{1.6}{0.6}{0}
    
    % Position Vector %
    \draw[darkred!100] (\cx,\cy) -- (\cix,\ciy);
%    \draw (\cx+\lllx,\cy+\llly) node[anchor=east] {$\rI(t)$};
    
    % FoV Cone %
    \cone{\cx}{\cy}{\aa}{\Lcone}{\Acone}{\Dcone}{color=red!100,opacity=0.10,shading=axis,shading angle=90,left color=red!100,right color=red!0};
    \draw (\labx+0.85,\laby+0.2) node[text width=2.75cm,align=left,anchor=east] {Position-triggered \\ body-fixed \\ $\losmax$-cone};
    \draw (\cx+\lllx,\cy+\llly) to[out=135,in=-45] (\labx,\laby);
    
    % Axes %
    \isoaxes{\cx}{\cy}{90+\aa}{1.5}{12}{$\bframe$}{$\bvec{e}_1$}{$\bvec{e}_2$}{$\bvec{e}_3$}
    \isoaxes{\cix}{\ciy}{90}{1.0}{15}{$\iframe$}{$\bvec{e}_1$}{$\bvec{e}_2$}{$\bvec{e}_3$}
    
    % FoV Dimension %
    \draw[color=black!100,<-] (\cx+\fact*\lx,\cy+\fact*\ly) arc(-90+\aa:-90+\aa+\da:2*\fact*\LL) node[anchor=west] {$\lambda(t)$};
    \draw[color=black!100,<-] (\cx+\fact*\llx,\cy+\fact*\lly) arc(-90+\aa+\aaa:-90+\aa+\aaa-\da:2*\fact*\LL);
    
    % Landing Site Label %
    \draw (\cix,\ciy) to[out=-90,in=-180] (0.5,-2.1) node[text width=3.5cm,align=left,anchor=west] {Landing site};
\end{tikzpicture}
    	\caption{Field of View State-Triggered Constraint}
    	\label{fig:stc_ex_fov}
   	\end{subfigure}
    \caption{Example applications where STCs are used to (a) limit the angle of attack $\alpha(t)$ at large dynamic pressures, and (b) to impose a field of view constraint that limits the line of sight angle $\lambda(t)$ to the landing site.}
\end{figure}

We now present an example application whose formulation within a continuous optimization framework is enabled by the state-triggered constraints introduced in~\sref{sec2:stc_formal_def}-\sref{sec2:stc_lcp}. Consider the problem of limiting the aerodynamic loads on a vehicle during a powered descent maneuver. On ascent, aerodynamic loads are often limited by imposing what are known as \qalpha~limits, where $q$ refers to dynamic pressure, and $\aoa$ refers to angle of attack. These constraints are typically valid only for small angles of attack, where aerodynamic loads are relatively easy to model. Unlike an ascent trajectory, a powered descent maneuver can exhibit a wider range of angles of attack, possibly exceeding $90\dg$ in cases where ``hopping'' maneuvers are permitted. Measures must therefore be taken to ensure that the vehicle does not operate in flight regimes with large model uncertainty. Specifically, we propose a simplified \qalpha~limit that enforces an angle of attack constraint \textit{only} at high dynamic pressures. Formally, we express this constraint using the following STC
\begin{equation} \label{eq:stc_ex_aoa_orig}
	\twonorm{\vB(t)} > \Vaoa \;\Rightarrow\; -\dotprod{\bvec{e}_1}{\vB(t)} \geq \cos{\aoamax}\twonorm{\vB(t)},
\end{equation}
where $\Vaoa\in\real_{++}$ is a speed above which the angle of attack is limited to $\aoa(t) \in\intee{0}{\aoamax}$. This STC is illustrated in Figure~\ref{fig:stc_ex_aoa}. Noting Remark~\ref{rem:stc_eq_ineq}, we convert the STC in~\eqref{eq:stc_ex_aoa_orig} into the following cSTC
\begin{subequations} \label{eq:stc_ex_aoa_imp}
	\begin{align}
		\haoa{\vI(t)}{\qIB(t)} &\definedas -\min\Big(\gaoa{\vI(t)}\,,0 \Big) \cdot \faoa{\vI(t)}{\qIB(t)} \leq 0\,, \label{eq:stc_ex_aoa_imp_a} \\
		\gaoa{\vI(t)} &\definedas \Vaoa-\twonorm{\vI(t)}, \\
    	\faoa{\vI(t)}{\qIB(t)} &\definedas \cos{\aoamax}\twonorm{\vI(t)}+\dotprod{\bvec{e}_1}{\cIB\big(\qIB(t)\big)\vI(t)}.
	\end{align}
\end{subequations}
\begin{remark}
    A range-triggered field of view constraint imposed between a body-fixed downward-looking camera and the landing site may be obtained by replacing~$\vI(t)$, $\Vaoa\,$, and~$\aoamax$ in~\eqref{eq:stc_ex_aoa_imp} with~$\rI(t)$, $\Rfov\,$, and~$\losmax\,$, respectively (see~Figure~\ref{fig:stc_ex_fov}). Such a constraint limits the line of sight angle to~$\lambda(t)\in\intee{0}{\losmax}$ when the vehicle is at distances greater than~$\Rfov$ away from the landing site.
\end{remark}
To conclude this section, we briefly discuss the shortcomings of three non-combinatorial alternatives to the cSTC proposed in~\eqref{eq:stc_ex_aoa_imp}:

\oursec{Alternative 1} Consider a naive temporally-scheduled implementation, in which the problem is first solved without \eqref{eq:stc_ex_aoa_orig}, and whose solution is used to determine the set of times~$\mathcal{T}$ over which the trigger condition in \eqref{eq:stc_ex_aoa_orig} is satisfied. The problem is then solved a second time with the constraint condition enforced over~$t\in\mathcal{T}$. This new solution now satisfies the constraint condition over~$t\in\mathcal{T}$, but does not necessarily satisfy the trigger condition for all~$t\in\mathcal{T}$. In fact, this solution may satisfy the trigger condition for times where~$t\notin\mathcal{T}$, thus necessitating a redefinition of~$\mathcal{T}$. Unlike cSTCs, this approach does not convey the relationship between the trigger and constraint conditions to the optimization, thus allowing this situation to persist.

\oursec{Alternative 2} Next, consider an implementation that replaces the STC in~\eqref{eq:stc_ex_aoa_orig} with $-\dotprod{\ex}{\vB(t)}/\twonorm{\vB(t)} \geq \cos{f_{\aoa}(t)}$, where $f_{\aoa}(t)\definedas f_{\aoa}\big(\twonorm{\vB(t)}\big)$ is a nonlinear scalar valued function that relates the maximum allowable angle of attack to the speed of the vehicle. This approach has two disadvantages: (i) it is less intuitive since the relationship between~$\aoa$ and~$q$ is embedded in~$f_{\aoa}(t)$, and (ii) obtaining a satisfactory~$f_{\aoa}(t)$ with proper numerical scaling may be difficult.

%%%%%%%%%%%%%%%%%%%%%%%%%%%%%%%%%%%%%%%%%%%%%%%%%%%%%%%%%%%%%%%%%%%%%%%%%%%%%%%%%%%%%%%%

 \boxing{t!}{problem}{prob:ncvx}{16cm}{Non-Convex Optimal Control Problem}{
 	\begin{tabular}{llll}
 		\underline{Cost Function}: &&& \\
 		&& $ {\begin{aligned}
        		&\underset{\tc,\,\tb,\,\TB(t)}{\text{minimize}}\;\; -\m(\tf) \\[1ex]
                &\text{s.t.}\quad\tc\in\intee{0}{\tcmax}
            \end{aligned}}$
        & $ {\begin{aligned}
        		& \\[1.5ex]
                &\text{See}\; \eqref{eq:bcs_initial_fit}
            \end{aligned}} $ \\
 		\underline{Boundary Conditions}: &&& \\
 		&& $ {\begin{aligned}
 		                 \m(\tig) &= \mig       && & \qIB(\tf)     &= \qIBf \\ 
 				        \rI(\tig) &= \prig(\tc) && & \rI(\tf)      &= \rIf \\
 					    \vI(\tig) &= \pvig(\tc) && & \vI(\tf)      &= \vIf \\
 					\omegaB(\tig) &= \omegaBig  && & \omegaB(\tf)  &= \omegaBf
         	\end{aligned}} $ 
         & $ {\begin{aligned}
            	&\text{See}\;\eqref{eq:bcs_final}\text{-}\eqref{eq:bcs_fit_polynomials}
              \end{aligned}} $ \\
 		\underline{Dynamics}: &&& \\
 		&& $ {\begin{aligned}
 					\dot{m}(t) &= -\mdotalpha\twonorm{\TB(t)} - \mdotbeta \\
 					\rIdot(t) &= \vI(t) \\
 					\vIdot(t) &= \frac{1}{m(t)} \cBI(t) \big(\TB(t) + \aeroB(t)\big) + \gI \\
 					\qIBdot(t) &= \frac{1}{2}\OMEGA{\omegaB(t)}\qIB(t) \\
                    \inertia \omegaBdot(t) &= \crossprod{\rTB}{\TB(t)} + \crossprod{\rCPB}{\aeroB(t)}- \crossprod{\omegaB(t)}{\inertia\omegaB(t)}
 				\end{aligned}} $ 
         & $ {\begin{aligned}
         		&\text{See}\; \eqref{eq:mass_dyn} \\
             	&\text{See}\; \eqref{eq:trans_dynamics_a} \\[1.5ex]
             	&\text{See}\; \eqref{eq:trans_dynamics_b} \\[1.5ex]
             	&\text{See}\; \eqref{eq:att_dyn_a} \\[1ex]
             	&\text{See}\; \eqref{eq:att_dyn_b}
             \end{aligned}} $ \\
        $ {\begin{aligned}
        		&\text{\underline{State Constraints}:} \\
                & \\
                & \\
                & \\
                & \\
		 		&\text{\underline{Control Constraints}:} \\
                & \\
                & \\
 				&\text{\underline{State-Triggered Constraints}:} \\
				& \\
             \end{aligned}} $
        && $ {\begin{aligned}
         		& \\
                \mdry &\leq m(t) \\
 				\tan{\glideslope}\twonorm{\Hgs\rI(t)} &\leq \dotprod{\ex}{\rI(t)} \\
 				\cos{\tiltmax} &\leq 1-2\twonorm{\Htilt\qIB(t)} \\
 				\twonorm{\omegaB(t)} &\leq \omegamax \\
 				& \\
                0 < \Tmin \leq \twonorm{\TB(t)} &\leq \Tmax \\
				\cos{\gimbalmax}\twonorm{\TB(t)} &\leq \dotprod{\ez}{\TB(t)} \\
                & \\
	         	\haoa{\vI(t)}{\qIB(t)} &\leq 0
             \end{aligned}} $
        & $ {\begin{aligned}
        		& \\
         		&\text{See}\; \eqref{eq:minmass} \\
             	&\text{See}\; \eqref{eq:glide_slope} \\
             	&\text{See}\; \eqref{eq:max_tilt} \\
             	&\text{See}\; \eqref{eq:wmax} \\
                & \\
         		&\text{See}\; \eqref{eq:thrust_bounds} \\
                &\text{See}\; \eqref{eq:max_gimbal} \\
                & \\
                &\text{See}\; \eqref{eq:stc_ex_aoa_imp} \\
             \end{aligned}} $ \\
         \null
 	\end{tabular}
}

\oursec{Alternative 3} Lastly, consider an implementation using two phases: the first with an angle of attack constraint and no velocity constraint, and the second with a velocity constraint and no angle of attack constraint. Further, assume that the terminal condition of the first phase is equated to the initial condition of the second, and that both phases are solved simultaneously. Such a multi-phase optimization approach ensures satisfaction of~\eqref{eq:stc_ex_aoa_orig}, and is well suited for applications where the temporal ordering and quantity of phases are known a priori. However, since the formulation presupposes a specific phase structure, this approach is not capable of introducing additional phases. In contrast, cSTCs can do so in order to achieve feasibility or improve optimality.

\subsection{Non-Convex Problem Statement} \label{sec2:nonconvex}

We now summarize the problem developed throughout this section. We assume a minimum-fuel objective function, but note that other objective functions can be readily used (e.g. a minimum-time problem would minimize $\tb$). The non-convex generalized powered descent guidance problem is summarized in Problem~\ref{prob:ncvx}.

%%%%%%%%%%%%%%%%%%%%%%%%%%%%%%%%%%%%%%%%%%%%%%%%%%%%%%%%%%%%%%%%%%%%%%%%%%%%%%%%%%%%%%%
%%%%%%%%%%%%%%%%%%%%%%%%%%%%%%%%%%%%%%%%%%%%%%%%%%%%%%%%%%%%%%%%%%%%%%%%%%%%%%%%%%%%%%%
%%%%%%%%%%%%%%%%%%%%%%%%%%%%%%%%%%%%%%%%%%%%%%%%%%%%%%%%%%%%%%%%%%%%%%%%%%%%%%%%%%%%%%%

\section{Convex Formulation} \label{sec:convex_formulation}

The successive convexification algorithm described in this section is designed to solve Problem~\ref{prob:ncvx} such that the converged solution (i) \textit{exactly} satisfies the original nonlinear dynamics, (ii) \textit{approximates} the state and control constraints by enforcing them only at a finite number of temporal nodes, and (iii) \textit{conservatively approximates} the optimality and feasibility of the problem by using a finite-dimensional representation of the infinite-dimensional control signal. The proposed algorithm works by iteratively solving a sequence of subproblems until a converged solution is attained. Each iteration consists of two steps: a \textit{propagation step} responsible for obtaining a subproblem, followed by a \textit{solve step} responsible for solving said subproblem to full optimality. Each subproblem is a convex approximation of Problem~\ref{prob:ncvx}, and each solve step results in a state and control trajectory, or \textit{iterate}. The propagation step in the first iteration is computed using a user-defined initialization trajectory (see~\sref{sec3:init}), whereas subsequent iterations perform said approximation about the trajectory obtained by the previous iteration's solve step. Since the solve step is executed using well understood algorithms (e.g. IPMs~\cite{nocedal2006numerical}), this section primarily focuses on the propagation step.

This section is organized as follows. In~\sref{sec3:gen_subproblem}, we outline a procedure to convert a \textit{free-final-time nonlinear continuous-time} optimal control problem into a general implementable \textit{fixed-final-time linear-time-varying discrete-time} convex parameter optimization subproblem~\cite{szmuk2018successive}. Specifically,~\sref{sec3:norm}-\sref{sec3:disc} discuss three analytical steps that comprise the propagation step, and~\sref{sec3:vctrl} introduces the virtual control and trust region modifications used to aid convergence. In~\sref{sec3:spec_subproblem}, we specialize the general subproblem to Problem~\ref{prob:ncvx}. Lastly, in~\sref{sec3:scvx_alg}, we describe two initialization strategies, and conclude by summarizing the proposed algorithm.

%%%%%%%%%%%%%%%%%%%%%%%%%%%%%%%%%%%%%%%%%%%%%%%%%%%%%%%%%%%%%%%%%%%%%%%%%%%%%%%%%%%%%%%

\subsection{General Implementable Convex Subproblem} \label{sec3:gen_subproblem}

This section assumes the following general optimal control problem
\begin{subequations} \label{eq:gen_ocp}
  \begin{align}
    &\underset{\tc,\,\tb,\,\uu(t)}{\text{minimize}} \;\; \Jxu{\tb}{\xx(\tf)}, \hspace{-2.2cm}&        &                           \label{eq:gen_ocp_a} \\
    \text{s.t.} \quad &\dot{\xx}(t) = \fxu{\xx(t)}{\uu(t)}\,,                 \hspace{-2.2cm}& \forall&\,t\in\intee{\tig}{\tf}\,, \label{eq:gen_ocp_b} \\
    &\hi{\zz(t)}{i} = 0\,,                                                    \hspace{-2.2cm}& \forall&\,t\in\intee{\tig}{\tf}\,,\; \forall\,i\in\setcvx\cup\setncvx\cup\setcstc, \label{eq:gen_ocp_c}
  \end{align}
\end{subequations}
where~$\xx(t)\in\real^{\ns}$ and~$\uu(t)\in\real^{\nc}$ denote respectively the state and control vectors,~$\zz(t)\definedas \big[\tc\;\, \tb\;\, \xx^T(t) \;\, \uu^T(t)\big]^T\in\real^{\nz}$ (see~\sref{sec2:notation} for definitions of~$\tc$ and~$\tb$),~$\fcn{\fcnJ}{\real_{++}\times\real^{\ns}}{\real}$ is a Mayer objective function~\cite{berkovitz_opt_75} that is convex in its arguments,~$\fcn{\fcnf}{\real^{\ns}\times\real^{\nc}}{\real^{\ns}}$ represents the continuous-time nonlinear dynamics, and~$\fcn{\fcncons_i}{\real^{\nz}}{\real}$ represent equality constraints imposed on the trajectory. The equality constraints in~\eqref{eq:gen_ocp_c} are assumed to be convex for~$i\in\setcvx\definedas\{1,\,\ldots\,,\np\}$, and non-convex for~$i\in\setncvx\definedas\{\np+1,\,\ldots\,,\np+\nq\}$ and~$i\in\setcstc\definedas\{\np+\nq+1,\,\ldots\,,\np+\nq+\nr\}$. The functions~$\fcnJ$ and~$\fcnf$ are assumed to be continuously differentiable, whereas each~$\fcncons_i$ is assumed to be at least once differentiable almost everywhere.

%%%%%%%%%%%%%%%%%%%%%%%%%%%%%%%%%%%%%%%%%%%%%%%%%%%%%%%%%%%%%%%%%%%%%%%%%%%%%%

\subsubsection{Normalization} \label{sec3:norm}

The normalization step converts the \textit{free-final-time} nonlinear continuous-time optimal control problem in~\eqref{eq:gen_ocp} into an equivalent \textit{fixed-final-time} nonlinear continuous-time problem. This is achieved by temporally normalizing the burn phase from~$t\in\intee{\tig}{\tf}$ to~$\tau\in\intee{0}{1}$, where~$\tnorm$ is the normalized burn phase time.
Using the chain rule, the nonlinear dynamics in~\eqref{eq:gen_ocp_b} can be rewritten as 
\begin{equation} \label{eq:tnorm_dyn}
	\xx'(t) \definedas \dd{}{\tnorm}\xx(t) = \dd{t}{\tnorm}\dd{}{t}\xx(t) = \left(\dd{t}{\tnorm}\right)\dot{\xx}(t).
\end{equation}
Defining the \textit{dilation factor}~$\dilation\definedas dt / d\tnorm\in\real_{++}$ 
and replacing~$\dot{\xx}(t)$ with the right-hand side of \eqref{eq:gen_ocp_b}, the temporally-normalized dynamics are expressed as
\begin{equation} \label{eq:normalized_dyn}
	\xx'(\tnorm) = \fxus{\xx(\tnorm)}{\uu(\tnorm)}{\dilation} \definedas \dilation\cdot\fxu{\xx(\tnorm)}{\uu(\tnorm)}.
\end{equation}
Since~$\tnorm\in\intee{0}{1}$, it follows that~$\dilation = \tb$. Thus, the temporal normalization of~\eqref{eq:gen_ocp_a} and~\eqref{eq:gen_ocp_c} is achieved by replacing the first argument of the cost function with~$\dilation$, and all subsequent~$\tf$ and~$t$ arguments with~$1$ and~$\tnorm$, respectively.

%%%%%%%%%%%%%%%%%%%%%%%%%%%%%%%%%%%%%%%%%%%%%%%%%%%%%%%%%%%%%%%%%%%%%%%%%%%%%%

\subsubsection{Linearization} \label{sec3:lin}

The linearization step converts the fixed-final-time \textit{nonlinear} continuous-time problem obtained in~\sref{sec3:norm} into a fixed-final-time \textit{linear-time-varying} continuous-time problem. By approximating non-convexities to first-order, the linearization step guarantees convexity of the subproblem.

\oursec{Dynamics} The right-hand side of~\eqref{eq:normalized_dyn} is approximated by a first-order Taylor series, evaluated about a reference trajectory denoted by~$\zzo(\tnorm)\definedas\big[ \tco\;\, \sso\;\, \xxo^T(\tnorm)\;\, \uuo^T(\tnorm)\big]^T$ for all~$\tnorm\in\intee{0}{1}$. 
The resulting linear-time-varying dynamics are given by
\begin{subequations} \label{eq:dynamics_lin}
	\begin{align}
    	\xx'(\tnorm) &\approx \Ac(\tnorm)\xx(\tnorm) + \Bc(\tnorm)\uu(\tnorm) + \Sc(\tnorm)\dilation + \rr(\tnorm), \\[2ex]
    	\Ac(\tnorm) &\definedas \dfsdx\,,\\
   		\Bc(\tnorm) &\definedas \dfsdu\,,\\
	    \Sc(\tnorm) &\definedas \dfsds\,,\\
	    \rr(\tnorm) &\definedas - \Ac(\tnorm)\xxo(\tnorm) - \Bc(\tnorm)\uuo(\tnorm).
	\end{align}
\end{subequations}

\oursec{Non-Convex State and Control Constraints} Using the assumption that the functions~$\stch_{i}(\cdot)$ are at least once differentiable almost everywhere, we define $\delta\zz(\tau)\definedas\zz(\tau)-\zzo(\tau)$ and approximate the constraints for each $i\in\setncvx$ in~\eqref{eq:gen_ocp_c} using a first-order Taylor series:
\begin{align}
      \hi{\zz(\tnorm)}{i} \approx \hi{\zzo(\tnorm)}{i} &+ \dhidz \delta\zz(\tnorm)\,.\quad \label{eq:constraint_ncvx_lin}
\end{align}

\oursec{Continuous State-Triggered Constraints} Similarly, since the trigger and constraint functions of each cSTC are assumed to be piecewise continuously differentiable, it follows that each~$\stch_{i}$ is at least once differentiable almost everywhere for all~$i\in\setcstc$.
However, due to the~$\min(\cdot)$ function in~\eqref{eq:stc_lcp}, the partial~$\partial\stch_{i}/\partial\zz$ is not well defined when~$\stcg_{i}(\cdot)=0$. To ensure that the constraint condition is not enforced when~$\stcg_{i}\big(\zzo(\tnorm)\big) = 0$, we define~$\partial\stch_{i}/\partial\zz$ to hold the same value as when~$\stcg_{i}\big(\zzo(\tnorm)\big) > 0$. Thus, the approximation is given as follows:
\begin{align}
	\hi{\zz(\tnorm)}{i} \approx
	\left\{\begin{aligned}
		\hi{\zzo(\tnorm)}{i}+\dhidz \delta\zz(\tnorm) &\,,\quad\text{if}\;\,\stcg_{i}\big(\zzo(\tnorm)\big) < 0\,, \\
                                                   0 &\,,\quad\text{otherwise.} \\
	\end{aligned}\right. \label{eq:constraint_cstc_lin}
\end{align}

%%%%%%%%%%%%%%%%%%%%%%%%%%%%%%%%%%%%%%%%%%%%%%%%%%%%%%%%%%%%%%%%%%%%%%%%%%%%%%

\subsubsection{Discretization} \label{sec3:disc}

The discretization step converts the fixed-final-time linear-time-varying \textit{continuous-time} problem obtained in~\sref{sec3:lin} into a fixed-final-time linear-time-varying \textit{discrete-time} parameter optimization problem. This step is critical in ensuring that the converged solution \textit{exactly} adheres to the prescribed \textit{nonlinear} dynamics. We begin by introducing~$\KK$ evenly spaced temporal nodes that divide the burn phase into~$\KK-1$ subintervals, and define the sets~$\setK\definedas\{1,\,2,\,\ldots\,,\,\KK\}$ and~$\setKm\definedas\{1,\,2,\,\ldots\,,\,\KK-1\}$. Each temporal node is associated with an index~$k\in\setK$, and corresponding normalized time~$\tauk = (k-1)/(\KK-1)$.

To proceed, the control signal must be projected from the infinite-dimensional space it inhabits to a finite-dimensional space suitable for numerical optimization. This can be done in numerous ways, including zero-order-hold (ZOH) and first-order-hold (FOH) interpolation. Alternatively, the state may be projected to a finite-dimensional space directly using pseudospectral methods~\cite{betts1998,DM_fahroo}. We have found that, compared to pseudospectral methods, ZOH and FOH interpolation yield sparsity patterns that noticeably decrease solve time. Our approach utilizes FOH interpolation because it (i) provides a noticeable increase in optimality when compared ZOH interpolation, and (ii) ensures that when the discrete-time control variables satisfy convex control constraints, the interpolated values follow suit. 
Formally, FOH interpolation represents the control signal over each subinterval~$k\in\setKm$ as
\begin{subequations} \label{eq:foh_def}
	\begin{gather}
	    \uu(\tnorm) = \lambdalk{\tnorm}\uuk+\lambdark{\tnorm}\uukp\,,\quad\forall\tnorm\in\intee{\tauk}{\taukp}\,,\\[2ex]
	\lambdalk{\tnorm}\definedas\frac{\taukp-\tnorm}{\taukp-\tauk}\,,\quad\lambdark{\tnorm}\definedas\frac{\tnorm-\tauk}{\taukp-\tauk},
	\end{gather}
\end{subequations}
where~$\uuk\definedas\uu(\tauk)$. Substituting~\eqref{eq:foh_def} into \eqref{eq:dynamics_lin}, we 
obtain the following for each subinterval~$k\in\setKm$:
\begin{equation}
	\xx'(\tnorm) = \Ac(\tnorm)\xx(\tnorm) + \lambdalk{\tnorm}\Bc(\tnorm)\uuk + \lambdark{\tnorm}\Bc(\tnorm)\uukp + \Sc(\tnorm)\dilation + \rr(\tnorm)\,,\quad\forall\,\tnorm\in\intee{\tauk}{\taukp}. \label{eq:norm_lti_dyn}
\end{equation}
The state transition matrix~$\stm{\xi}{\tauk}$ for~$\xi\in\intee{\tauk}{\taukp}$ associated with~\eqref{eq:norm_lti_dyn} is given by
\begin{equation} \label{eq:stm_dyn}
	\stm{\xi}{\tauk} =  I_{\ns\times\,\ns} + \int_{\tauk}^{\xi}{A(\zeta)\,\stm{\zeta}{\tauk}d\zeta}\,.
\end{equation}
Denoting the discrete-time state vectors by~$\xxk\definedas\xx(\tauk)$, the inverse and transitive properties of~$\stm{\cdot}{\cdot}$~\cite{Antsaklis2007} are used to obtain the following discrete-time solution to~\eqref{eq:norm_lti_dyn} for each~$k\in\setKm$:
\begin{subequations} \label{eq:ltv_disc}
	\begin{align}
		\xxkp &= \Ad\xxk+\Bdm\uuk+\Bdp\uukp+\Sd\dilation+\rd, \label{eq:ltv_disc_a} \\[2ex]
      	\Ad  &\definedas \stm{\taukp}{\tauk}, \label{eq:ltv_disc_b} \\
      	\Bdm &\definedas \Ad\int_{\tauk}^{\taukp}{\stminv{\xi}{\tauk}\Bc(\xi)\lambdalk{\xi}d\xi}, \label{eq:ltv_disc_c} \\
      	\Bdp &\definedas \Ad\int_{\tauk}^{\taukp}{\stminv{\xi}{\tauk}\Bc(\xi)\lambdark{\xi}d\xi}, \label{eq:ltv_disc_d} \\
      	\Sd  &\definedas \Ad\int_{\tauk}^{\taukp}{\stminv{\xi}{\tauk}\Sc(\xi)d\xi}, \label{eq:ltv_disc_e} \\
      	\rd  &\definedas \Ad\int_{\tauk}^{\taukp}{\stminv{\xi}{\tauk}\rr(\xi)d\xi}. \label{eq:ltv_disc_f}
	\end{align}
\end{subequations}

In implementation, the previous iteration's solve step generates~$\tco$, $\sso$, $\xxko\definedas\xxo(\tauk)$, and~$\uuko\definedas\uuo(\tauk)$ for all~$k\in\setK$. The~$\uuko$'s and~\eqref{eq:foh_def} are used to obtain $\uuo(\tnorm)$ over~$\tau\in\intee{0}{1}$, and~\eqref{eq:normalized_dyn}, \eqref{eq:ltv_disc_b}, and the integrands in~\eqref{eq:ltv_disc_c}-\eqref{eq:ltv_disc_f} are computed simultaneously for each~$k\in\setKm$ using the intermediate quantity~$\stm{\xi}{\tauk}$. When the propagation reaches $\taukp$, the quantity in~\eqref{eq:ltv_disc_b} can be left multiplied against the integrands of~\eqref{eq:ltv_disc_c}-\eqref{eq:ltv_disc_f} to obtain the final values of~$\Ad$, $\Bdm$, $\Bdp$, $\Sd$, and~$\rd$. 

The integration of~\eqref{eq:normalized_dyn} is initialized with~$\xxko$, and is analogous to a multiple shooting method. We have observed that this multiple shooting strategy improves the convergence behavior of the algorithm by keeping $\zzo(\tnorm)$ closer to the path obtained by the constrained optimization problem. In contrast, since the dynamics are nonlinear and may be unstable, a single shooting method is more susceptible to poor initializations, since the dynamics have more time to evolve away from feasibility. Thus, the propagation and solve steps are designed to play complementary roles, whereby the former relates the discrete-time optimization problem to the continuous-time physics in the absence of constraints, and the latter helps reset the shooting method at more frequent temporal intervals while taking into account the constraints. 

\begin{remark}
	The initial condition of the~$\ith{k}{th}$ subinterval depends only on~$\zzko$ obtained by the solve step of the previous iteration. It \underline{does not} depend on the propagation of the~$\ith{(k-1)}{th}$ subinterval. Therefore, the propagation step can be parallelized such that all subintervals~$k\in\setKm$ are computed simultaneously.
\end{remark}

The state and control constraints are enforced at the temporal nodes~$\tauk$ for all~$k\in\setK$. This discretization choice is adopted for simplicity, and does not guarantee the absence of inter-node state constraint and non-convex control constraint violations.

We conclude by noting that the number of temporal nodes and interpolation scheme do not affect the \textit{accuracy} of the converged solution with respect to the nonlinear dynamics. Instead, these choices affect the \textit{optimality} and, in extreme cases, the \textit{feasibility} of the solution. Simply put, a coarser or less expressive interpolation results in a finite-dimensional control signal with fewer degrees of freedom. This yields a problem that is implicitly more constrained, and thus a solution that is generally less optimal.

%%%%%%%%%%%%%%%%%%%%%%%%%%%%%%%%%%%%%%%%%%%%%%%%%%%%%%%%%%%%%%%%%%%%%%%%%%%%%%

\subsubsection{Virtual Control \& Trust Region Modifications} \label{sec3:vctrl}

The process outlined in~\sref{sec3:norm}-\sref{sec3:disc} results in an implementable convex parameter optimization subproblem that may suffer from \textit{artificial infeasibility} and \textit{artificial unboundedness}. These issues are addressed by the virtual control and trust region modifications, respectively. 

\oursec{Artificial Infeasibility \& Virtual Control} Consider linearizing and discretizing Problem~\ref{prob:ncvx} about a burn time~$\sso=0$ (or some small value). From~\eqref{eq:dynamics_lin}, \eqref{eq:stm_dyn}, and~\eqref{eq:ltv_disc}, it follows that~\eqref{eq:ltv_disc_a} does not admit a solution for arbitrary~$\xxkp\neq\xxk$. This remains true even if Problem~\ref{prob:ncvx} admits a feasible solution, resulting in a condition we term \textit{artificial infeasibility.}
To mitigate artificial infeasibility, we augment~\eqref{eq:ltv_disc_a} with a virtual control term~$\nuk\in\real^{\ns}$ for all~$k\in\setKm$, 
\begin{equation} \label{eq:ltv_disc_relaxed}
	\xxkp = \Ad\xxk+\Bdm\uuk+\Bdp\uukp+\Sd\dilation+\rd+\nuk\,,
\end{equation}
and add the following penalty term to~\eqref{eq:gen_ocp_a},
\begin{equation} \label{eq:J_vc}
	\Jvc{\nubig} \definedas \wvc\sum_{k\in\setKm}{\onenorm{\nuk}}\,,
\end{equation}
where~$\wvc\in\real_{++}$ is a large weight. The virtual control modification guarantees that each subproblem has a non-empty feasible set, and thus ensures that the convergence process is not obstructed. A 1-norm minimization is employed in~\eqref{eq:J_vc} to encourage sparsity in the vectors~$\nuk$. Upon successful convergence, the virtual control terms are zero, and~\eqref{eq:ltv_disc_relaxed} is equivalent to~\eqref{eq:ltv_disc}.

\oursec{Artificial Unboundedness \& Trust Region} The second issue that may arise is that of \textit{artificial unboundedness}, which occurs when the \textit{linearized} constraints permit the cost of a subproblem to be minimized indefinitely. This issue is mitigated by adding the following soft quadratic trust region to the cost function,
\begin{equation} \label{eq:J_tr}
	\Jtr{\zzo}{\zz} \definedas \sum_{k\in\setK}{\delta\zzk^T\,\Wtr\,\delta\zzk}\,,
\end{equation}
where~$\delta\zzk\definedas\delta\zz(\tauk)$ for all~$k\in\setK$, and~$\Wtr\in\spd{\nz}$ is a symmetric positive definite weighting matrix. The trust region modification also serves to ensure that the solve step does not venture excessively far from the reference trajectory used in the propagation step. 

%%%%%%%%%%%%%%%%%%%%%%%%%%%%%%%%%%%%%%%%%%%%%%%%%%%%%%%%%%%%%%%%%%%%%%%%%%%%%%

\subsection{Specialized Implementable Convex Subproblem} \label{sec3:spec_subproblem}

To specialize the subproblem developed in~\sref{sec3:gen_subproblem} to Problem~\ref{prob:ncvx}, we define the objective function~\eqref{eq:gen_ocp_a} to be the minimum-fuel objective~$\Jxu{\dilation}{\xxkf} \definedas -\mkf$. We define~$\xx(t) \definedas \big[\m(t)\;\, \rI^T(t)\;\, \vI^T(t)\;\, \qIB^T(t)\;\, \omegaB^T(t)\big]^T$ and~$\uu(t) \definedas \TB(t)$, and concatenate~\eqref{eq:mass_dyn}-\eqref{eq:att_dyn} to form the corresponding dynamics.
In accordance with~\eqref{eq:constraint_ncvx_lin}, the thrust lower bound constraint in~\eqref{eq:thrust_bounds} is approximated using the following form for each~$k\in\setK$
\begin{equation} \label{eq:lin_thrust_lower_bound}
	\ftlbk + \Htlbk\,\delta\zzk \leq 0.
\end{equation}
Similarly, in accordance with~\eqref{eq:constraint_cstc_lin}, the cSTC in~\eqref{eq:stc_ex_aoa_imp} is approximated using the following form for each~$k\in\setK$
\begin{equation} \label{eq:lin_aoa_cstc}
	\haoak + \Haoak\,\delta\zzk \leq 0.
\end{equation}

Problem~\ref{prob:cvx} presents a summary of the specialized convex parameter optimization subproblem used in our algorithm. This problem is primarily concerned with solving for the states~$\xxk$, controls~$\uuk$, and time dilation~$\dilation$ associated with the \textit{burn phase} of the trajectory. The only aspect of the \textit{coast phase} optimized by Problem~\ref{prob:cvx} is the coast time~$\tc$, which in turn determines the ignition time position~$\rIki$ and velocity~$\vIki$ via~\eqref{eq:bcs_initial_fit} and~\eqref{eq:bcs_fit_polynomials}.

\boxing{ht!}{problem}{prob:cvx}{16cm}{Convex Parameter Optimization Subproblem}{
	\begin{tabular}{llll}
		\underline{Cost Function}:       &&& \\
		                                 & ${\begin{aligned}
                                         	 &\underset{\tc,\,\dilation,\,\uuk,\,\nuk}{\text{minimize}}\; -\mk+\Jvc{\nubig}+\Jtr{\zzo}{\zz} \\[1ex]
                                             &\text{s.t.}\quad\tc\in\intee{0}{\tcmax}
                                            \end{aligned}}$
                                         &
                                         & $ {\begin{aligned}
                                         	\text{See}\; &\eqref{eq:J_vc}\,\&\,\eqref{eq:J_tr} \\[0.7ex]
                                            \text{See}\; &\eqref{eq:bcs_initial_fit}
                                         \end{aligned}} $ \\
		\underline{Boundary Conditions}: & & &\\
		                                 & $ {\begin{aligned}
				                                \mki      &= \mig       && & \qIBkf        &= \qIBf \\ 
					                            \rIki     &= \prig(\tc) && & \rIkf         &= \rIf \\
					                            \vIki     &= \pvig(\tc) && & \vIkf         &= \vIf \\
					                            \omegaBki &= \omegaBig  && & \omegaBkf     &= \omegaBf
				                              \end{aligned}} $
                                         & 
                                         & $ {\begin{aligned}
                                             & \text{See}\;\eqref{eq:bcs_initial}\text{-}\eqref{eq:bcs_fit_polynomials}
                                             \end{aligned}} $ \\ 
		$ {\begin{aligned}
        	& \text{\underline{Dynamics}:} \\
            & \\
            &
           \end{aligned}} $
        & $ {\begin{aligned}
            & \\
            & \xxkp = \Ad\xxk+\Bdm\uuk+\Bdp\uukp+\Sd\dilation+\rd+\nuk \\
            &
           \end{aligned}}$
        & $ {\begin{aligned}
            & \\
            & \forall k\in\setKm \\
            &
           \end{aligned}} $
        & $ {\begin{aligned}
            & \\
            \text{See}\, &\eqref{eq:normalized_dyn}\text{-}\eqref{eq:dynamics_lin}, \\
                         & \eqref{eq:foh_def}\text{-}\eqref{eq:ltv_disc}%\,\& \\
                                            %& \eqref{eq:xu_def}\text{-}\eqref{eq:specialized_dyn}
           \end{aligned}} $ \\
		$ {\begin{aligned}
        	& \text{\underline{State Constraints}:} \\
            & \\
            & \\
            & \\
            & \\
            & \text{\underline{Control Constraints}:} \\
            & \\
            & \\
            & \\
           	& \text{\underline{State-Triggered Constraints}:} \\
            &
           \end{aligned}} $
		& $ {\begin{aligned}
        			& \\
					\mdry &\leq \mk \\
					\tan{\glideslope}\twonorm{\Hgs\rIk} &\leq \dotprod{\ex}{\rIk} \\
					\cos{\tiltmax} &\leq 1-2\twonorm{\Htilt\qIBk} \\
					\twonorm{\omegaBk} &\leq \omegamax \\
                    & \\
                    \twonorm{\uuk} &\leq \Tmax \\
					\cos{\gimbalmax}\twonorm{\uuk} &\leq \dotprod{\ez}{\uuk} \\
					\ftlbk+\Htlbk\,\delta\zzk &\leq 0 \\
                    & \\
					\haoak+\Haoak\,\delta\zzk &\leq 0
			\end{aligned}} $
        & $ {\begin{aligned}
        		& \\
                & \\[2ex]
            	&\forall k\in\setK \\[1ex]
                & \\
                & \\
                & \\[0.6ex]
            	&\forall k\in\setK \\
                & \\
                & \\
                &\forall k\in\setK
        	\end{aligned}} $
        & $ {\begin{aligned}
        		& \\
	            \text{See}\; &\eqref{eq:minmass} \\
    	        \text{See}\; &\eqref{eq:glide_slope} \\
                \text{See}\; &\eqref{eq:max_tilt} \\
            	\text{See}\; &\eqref{eq:wmax} \\
                & \\
            	\text{See}\; &\eqref{eq:thrust_bounds} \\
            	\text{See}\; &\eqref{eq:max_gimbal} \\
            	\text{See}\; &\eqref{eq:thrust_bounds}\,\&\,\eqref{eq:lin_thrust_lower_bound} \\
                & \\
				\text{See}\; &\eqref{eq:stc_ex_aoa_imp}\,\&\,\eqref{eq:lin_aoa_cstc}
            \end{aligned}} $ \\
            \null
	\end{tabular}
  }

%%%%%%%%%%%%%%%%%%%%%%%%%%%%%%%%%%%%%%%%%%%%%%%%%%%%%%%%%%%%%%%%%%%%%%%%%%%%%%

\subsection{Successive Convexification Algorithm} \label{sec3:scvx_alg}

\subsubsection{Initialization} \label{sec3:init}

We consider two initialization approaches: \textit{straight-line initialization} and \textit{3-DoF initialization}. Both approaches assume~$\tco=0$ and a user-specified initial guess for~$\sso$. The former is equivalent to assuming~$\tin=\tig$, whereas the latter is equivalent to guessing the burn time~$\tb$ (see Figure~\ref{fig:traj_timeline}). The nature of the powered descent guidance problem is such that~$\sso$ can typically be guessed accurately as a function of distance to the landing site, initial velocity, and available thrust. In our experience, the proposed algorithm is able to handle a wide range of initialization values for~$\sso$, although poor guesses may lead to increased convergence time.

The straight-line initialization approach constructs the discrete-time state trajectory~$\xxko$ by linearly interpolating the state at each temporal node between the ignition and final states. The control trajectory~$\uuko$ is assumed to oppose the gravitational force at each temporal node. The initialization assumes an initial attitude~$\qidentity$ and final mass~$\mdry$, since these quantities are not known a priori. Formally, the state and control for each~$k\in\setK$ is computed as follows
\begin{subequations}
	\begin{align}
		\xxko = &\left(\dfrac{\KK-k}{\KK-1}\right)\xxoig+\left(\dfrac{k-1}{\KK-1}\right)\xxof, \\
        \uuko = -&\left(\dfrac{\KK-k}{\KK-1}\right)\mig\gI-\left(\dfrac{k-1}{\KK-1}\right)\mdry\gI,
	\end{align}
\end{subequations}
where~$\xxoig\definedas\big[\mig\;\,\rIig^T\;\,\vIig^T\;\,\qidentity^T\;\,\omegaBig^T\big]^T$, and~$\xxof\definedas\big[\mdry\;\,\rIf^T \vIf^T\;\,\qIBf^T\;\,\omegaBf^T\big]^T$. If a value other than~$\qidentity$ is assumed for the initial attitude, this approach can be improved by interpolating the quaternion states using Spherical Linear Interpolation (SLERP)~\cite{Hanson2006}.

The 3-DoF initialization approach constructs the state and control trajectories using a solution obtained from a convex 3-DoF guidance problem~\cite{behcetjgcd07}. The 3-DoF problem is solved using the same number of temporal nodes,~$\KK$, as in the 6-DoF problem. This problem may be solved once using a user-defined burn time, or in conjunction with a line-search that optimizes the burn time, and thus generates~$\sso$. The mass, position, and velocity components of~$\xxko$ and the controls~$\uuko$ are obtained directly from the 3-DoF solution. The attitude is computed such that the vertical axis of the vehicle is aligned with~$\uuko$, and the angular velocity is obtained by inverting~\eqref{eq:att_dyn_a}. 

Unsurprisingly, the initialization approach can have a significant impact on the converged solution attained by the algorithm. In our work, we have found that neither approach offers a clear and consistent advantage over the other, (see~\sref{sec4:timing}). We ultimately regard the initialization approach as a design choice.

%%%%%%%%%%%%%%%%%%%%%%%%%%%%%%%%%%%%%%%%%%%%%%%%%%%%%%%%%%%%%%%%%%%%%%%%%%%%%%
 
\subsubsection{Algorithm} \label{sec3:alg}

 \begin{algorithm}[b]
     \caption{\hspace{-0.1cm}\textbf{.} 6-DoF Powered Descent Guidance Successive Convexification Algorithm}
     \label{alg:scvx}
     \begin{algorithmic}[1] % [1] means where line numbering starts
         \State Initialize $\{\tco,\sso,\xxo,\uuo\}$
             \While{not converged}
                 \State Compute $\{\Ad,\Bdp,\Bdm,\Sd,\rd\}\;\;\forall\,k\in\setKm$ according to~\sref{sec3:disc}\Comment{Propagation Step}
                 \State Solve Problem~\ref{prob:cvx} to obtain $\{\tc,\dilation,\xx,\uu,\nubig\}$ \Comment{Solve Step}
                 \If{$ \Jvc{\nubig} < \epsvc$ and $\Jtr{\zz}{\zzo} < \epstr$} \Comment{Convergence Criteria}
                 \State converged
                 \EndIf
                 \State $\{\tco,\sso,\xxo,\uuo\} \gets \{\tc,\dilation,\xx,\uu\}$
             \EndWhile\label{euclidendwhile}
         \State \textbf{return} $\{\tc,\dilation,\xx,\uu\}$
     \end{algorithmic}
 \end{algorithm}

The algorithm is initialized using one of the two approaches discussed in~\sref{sec3:init}. For each iteration, the algorithm performs a propagation step to compute~$\Ad$, $\Bdp$, $\Bdm$, $\Sd$, and~$\rd$ for all~$k\in\setKm$, followed by a solve step that solves the convex second-order cone programming subproblem summarized in Problem~\ref{prob:cvx}. The process terminates when~$\Jvc{\nubig}<\epsvc$ and~$\Jtr{\zz}{\zzo}<\epstr$, where~$\epsvc,\,\epstr\in\real_{++}$ are user-specified convergence tolerances. Satisfaction of the convergence criteria ensures that the attained solution eliminates the first order terms~$\delta\zzk$ generated by the approximation without using virtual control. If an iterate satisfies the convergence criteria, then feasibility of the corresponding subproblem implies that the iterate exactly satisfies the nonlinear dynamics of Problem~\ref{prob:ncvx} for all $t\in\intee{\tig}{\tf}$, and satisfies the state and control constraints of Problem~\ref{prob:ncvx} at each temporal node. However, failure to converge does not necessarily imply that Problem~\ref{prob:ncvx} is infeasible. A summary of the algorithm is provided in Algorithm~\ref{alg:scvx}.

We note that prior to convergence, the iterates may not be feasible with respect to Problem~\ref{prob:ncvx} due to the linearization used in the propagation step. This statement holds true even though each convex subproblem is designed to be feasible through the use of virtual control.

Finally, we highlight that no convergence guarantees are presented in this paper. However, we have found that Algorithm~\ref{alg:scvx} works well in practice, and note its similarity to~\cite{mao2016successive}, which does guarantee convergence to a local optima of the original problem \textit{when} the converged solution requires no virtual control.

%%%%%%%%%%%%%%%%%%%%%%%%%%%%%%%%%%%%%%%%%%%%%%%%%%%%%%%%%%%%%%%%%%%%%%%%%%%%%%%%%%%%%%%
%%%%%%%%%%%%%%%%%%%%%%%%%%%%%%%%%%%%%%%%%%%%%%%%%%%%%%%%%%%%%%%%%%%%%%%%%%%%%%%%%%%%%%%
%%%%%%%%%%%%%%%%%%%%%%%%%%%%%%%%%%%%%%%%%%%%%%%%%%%%%%%%%%%%%%%%%%%%%%%%%%%%%%%%%%%%%%%

\section{Numerical Results} \label{sec:numerical_results}

In this section, we present simulation results that demonstrate the proposed successive convexification algorithm, while highlighting the principal contributions of this paper. In~\sref{sec4:exI}, \sref{sec4:exII}, and~\sref{sec4:exIII} we present case studies that respectively illustrate the effects of the aerodynamic models introduced in~\sref{sec2:aero}, the state-triggered constraints introduced in~\sref{sec2:stc}, and the free-ignition-time modification introduced in~\sref{sec2:fit}. In~\sref{sec4:timing}, we provide performance and timing results. To present the results, we introduce the problem feature labels given in Table~\ref{tab:prob_names}.

The simulations are designed around a notional non-dimensionalized scenario with time, length, and mass units~$\TU$, $\LU$, and~$\MU$. Each scenario is defined by the problem parameters in Tables~\ref{tab:params_ex} and the initial position and velocity vectors defined in each subsequent subsection. For the sake of illustration, the initial conditions used in the case studies in~\sref{sec4:exI}-~\sref{sec4:exIII} define in-plane maneuvers, whereas those used to generate the timing results in~\sref{sec4:timing} define more computationally intensive out-of-plane maneuvers.

% see https://www.aiaa.org/figuretableguidelines/ %
\vspace{0.25cm}
\begin{table}[h!]
  \centering
  \caption{Generalized Powered Descent Guidance Problem Features}
  \begin{tabular}{c c p{13cm}}
  	\hhline{===}
    Feature & Section                & Description \\
    \hline
    (B)     & \sref{sec2:6dof}       & Baseline problem setup with no aerodynamics and straight-line initialization. \\
    (FI)    & \sref{sec2:fit}        & Includes the free-ignition-time modification shown in Figure~\ref{fig:traj_timeline}, where initial position and velocity are restricted to a free fall trajectory prior to engine ignition. \\
    (SA)    & \sref{sec2:aero}       & Includes the spherical aerodynamic model from Figure~\ref{fig:aerodynamic_model_a}. \\
    (EA)    & \sref{sec2:aero}       & Includes the ellipsoidal aerodynamic model from Figure~\ref{fig:aerodynamic_model_b}. \\
    (ST)    & \sref{sec2:stc_ex_app} & Includes the state-triggered constraint from Figure~\ref{fig:stc_ex_aoa}, and introduced in~\eqref{eq:stc_ex_aoa_orig}. \\
    (3I)    & \sref{sec3:init}       & Uses the 3-DoF initialization in lieu of straight-line initialization. Straight-line initialization is implied in the absence of this feature. \\
    \hhline{===}
    \label{tab:prob_names}
  \end{tabular}
\end{table}
\vspace{-0.25cm}

\vspace{0.25cm}
\newcommand{\tabspace}{\hspace{-0.35cm}}
\newcommand{\tabsep}{\null\hspace{0.5cm}\null}
\begin{table}[htb]
  \centering
  \caption{Problem Parameters}
  \begin{tabular}{lrllclrll}
  	\hhline{====~====}
    \textbf{Parameter}      &   & \tabspace\textbf{Value}                                 & \textbf{Units} &\tabsep&
    \textbf{Parameter}      &   & \tabspace\textbf{Value}                                 & \textbf{Units} \\
    \hhline{----~----}
    $\gI$                   &$-$& \tabspace$\ex$                                          & $\LU/\TU^2$ &\tabsep&
	$\Vaoa$                 &   & \tabspace$2.0$                                          & $\LU/\TU$ \\
    $\density$              &   & \tabspace$1.0$                                          & $\MU/\LU^3$ &\tabsep&
	$\aoamax$               &   & \tabspace$3.0$                                          & $\dg$ \\
    $\inertia$              &   & \tabspace$0.168\cdot\diag{\big[2\text{e-}2\;1\;1\big]}$ & $\MU\cdot\LU^2$ &\tabsep&
	$\mdry$                 &   & \tabspace$1.0$                                          & $\MU$ \\
    $\Pamb$                 &   & \tabspace$0.0$                                          & $\MU/\TU^2/\LU$ &\tabsep&
	$\mig$                  &   & \tabspace$2.0$                                          & $\MU$ \\
    $\Anoz$                 &   & \tabspace$0.0$                                          & $\LU^2$ &\tabsep&
	$\rIkf$                 &   & \tabspace$\bvec{0}_{3 \times 1}$                        & $\LU$ \\
    $\rCPB$                 &   & \tabspace$0.05\cdot\ex$                                 & $\LU$ &\tabsep&
	$\vIkf$                 &$-$& \tabspace$0.1\cdot\ex$                                  & $\LU/\TU$ \\
    $\rTB$                  &$-$& \tabspace$0.25\cdot\ex$                                 & $\LU$ &\tabsep&
	$\omegaBki,\,\omegaBkf$ &   & \tabspace$\bvec{0}_{3 \times 1}$                        & $\dg/\TU$ \\
    $\Isp$                  &   & \tabspace$30.0$                                         & $\TU$ &\tabsep&
	$\qIBkf$                &   & \tabspace$\qidentity$                                   & - \\
    $\tiltmax$              &   & \tabspace$90.0$                                         & $\dg$ &\tabsep&
	$\KK$                   &   & \tabspace$20$                                           & - \\
	$\omegamax$             &   & \tabspace$28.6$                                         & $\dg/\TU$ &\tabsep&
	$\wvc$                  &   & \tabspace$1\text{e+}4$                                  & - \\
    $\glideslope$           &   & \tabspace$75.0$                                         & $\dg$ &\tabsep&
    $\Wtr$                  &   & \tabspace$0.5$                                          & - \\
	$\gimbalmax$            &   & \tabspace$20.0$                                         & $\dg$ &\tabsep&
	$\epsvc$                &   & \tabspace$1\text{e-}8$                                  & - \\
    $\Tmin$                 &   & \tabspace$1.5$                                          & $\MU\cdot\LU/\TU^2$ &\tabsep&
    $\epstr$                &   & \tabspace$5\text{e-}4$                                  & - \\
    $\Tmax$                 &   & \tabspace$6.5$                                          & $\MU\cdot\LU/\TU^2$ &\tabsep&
    $\sso$                  &   & \tabspace$5.0$                                          & - \\
	\hhline{====~====}
    \label{tab:params_ex}
  \end{tabular}
\end{table}

%%%%%%%%%%%%%%%%%%%%%%%%%%%%%%%%%%%%%%%%%%%%%%%%%%%%%%%%%%%%%%%%%%%%%%%%%%%%%%%%%%

\subsection{Aerodynamic Models Case Study} \label{sec4:exI}

\begin{figure}[t!]
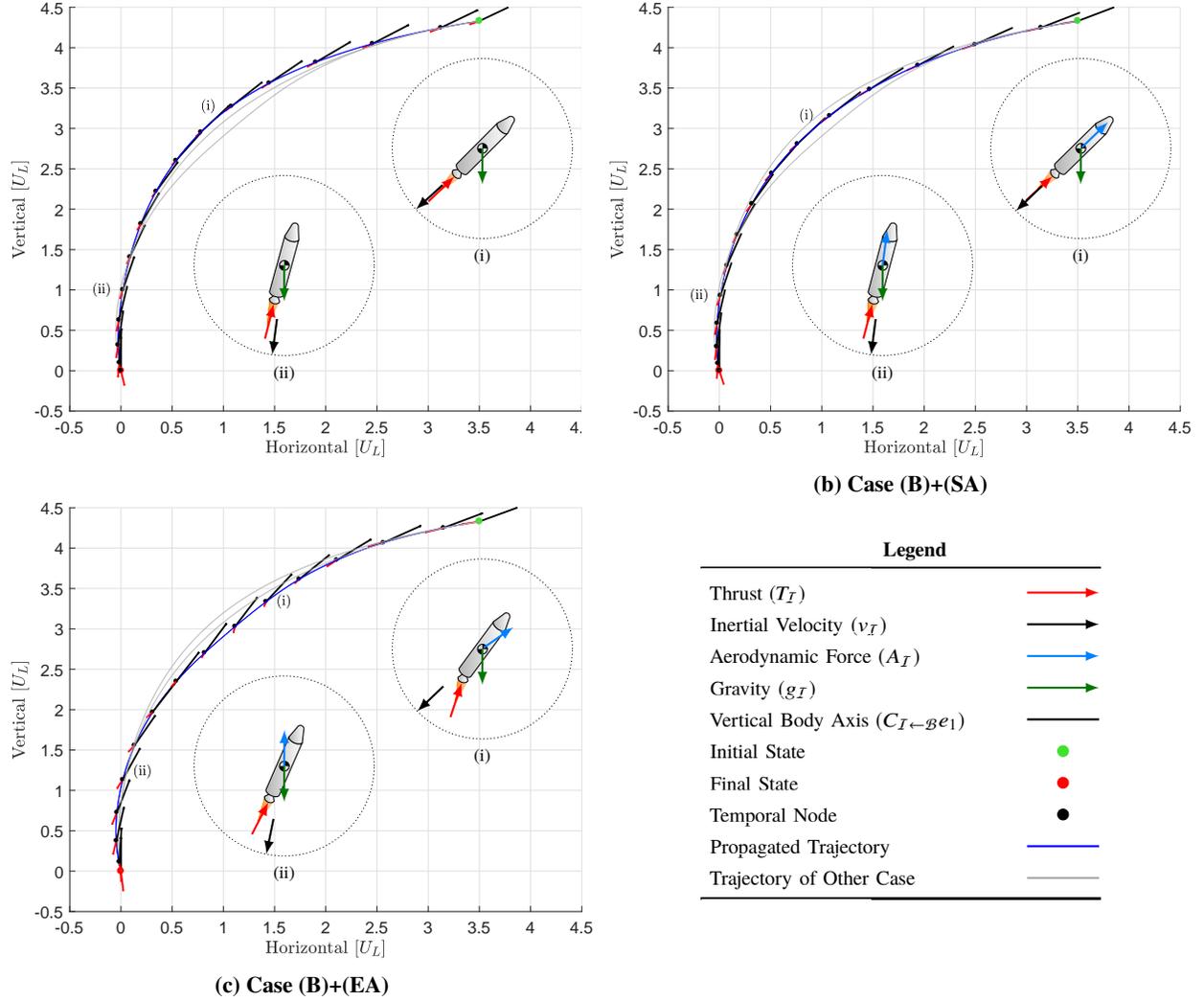

  \begin{subfigure}{0.5\textwidth}
      \centering
      \input{tikz_rocket_I.tex}
      \caption{Case~(B)}
      \label{fig:exI_nostc_nodrag}
  \end{subfigure} \hfil
  \begin{subfigure}{0.5\textwidth}
      \centering
      \input{tikz_rocket_II.tex}
      \caption{Case~(B)+(SA)}
      \label{fig:exI_nostc_sphere}
  \end{subfigure} \\
  \begin{subfigure}{0.5\textwidth}
  	  \centering
      \input{tikz_rocket_III.tex}
      \caption{Case~(B)+(EA)}
      \label{fig:exI_nostc_elip}
  \end{subfigure}  \hfil
  \begin{subfigure}{0.5\textwidth}
      \centering
      %\documentclass{article}
%\input{../commands.tex}
%\input{../mikz.tex}
%
%\begin{document}

\vspace{-0.5cm}
\newcommand{\legendspace}{-0.15cm}
\begin{tikzpicture}    
\tikzstyle{every node}=[font=\footnotesize]
\matrix[table] (mat11) 
{
\hspace{2.4cm} \textbf{Legend} & \\
%\textbf{Object} & \hspace{7.75em} \textbf{Style} \\
\hline
\text{Thrust } (\TI) & \draw[->,thick,color=ucol] (1.0,0.1) -- (2.0,0.1); \\[\legendspace]
\text{Inertial Velocity } (\vI) &   \draw[->,thick,color=vcol] (1.0,0.1) -- (2.0,0.1); \\[\legendspace]
\text{Aerodynamic Force } (\aeroI) &  \draw[->,thick,color=acol] (1.0,0.1) -- (2.0,0.1); \\[\legendspace]
\text{Gravity } (\gI) &  \draw[->,thick,color=gcol] (1.0,0.1) -- (2.0,0.1); \\[\legendspace]
\text{Vertical Body Axis } (\cBI\ex) &  \draw[thick,color=bcol] (1.0,0.1) -- (2.0,0.1); \\[\legendspace]
\text{Initial State} & \draw[fill,color=scol] (1.5,0.1) circle (0.075); \\[\legendspace]
\text{Final State} & \draw[fill,color=ucol] (1.5,0.1) circle (0.075); \\[\legendspace]
\text{Temporal Node} & \draw[fill] (1.5,0.1) circle (0.075);\\[\legendspace]
\text{Propagated Trajectory} & \draw[thick,color=blue!100] (1.0,0.1) -- (2.0,0.1); \\[\legendspace]
\text{Trajectory of Other Case} & \draw[thick,color=ocol] (1.0,0.1) -- (2.0,0.1);\\[0.2ex]
\hline \\
};
\end{tikzpicture}

%\end{document}
  \end{subfigure}
  \caption{Aerodynamic Models Case Study: Trajectories for the baseline (B), spherical aerodynamic model (B)+(SA) and ellipsoidal aerodynamic model (B)+(EA) cases with~$\tin=\tig$, $\rI(\tin)=[4.33 \ 3.5 \ 0.0 ]^T\;\LU$ and~$\vI(\tin) =[-0.5\ -2.5\ 0.0]^T\;\LU/\TU$. All force vectors are normalized to show direction only.}
  \label{fig:exI_trjs}
\end{figure}

In this case study, we solve three otherwise identical powered descent guidance problems, assuming no aerodynamic effects, a spherical aerodynamic model, and an ellipsoidal aerodynamic model. These problems are labeled using (B), (B)+(SA) and (B)+(EA), respectively. In each problem, the vehicle begins above and east of the landing pad, traveling west at a shallow flight path angle. These initial conditions are given by $\rI(\tin) = [4.33\ 3.5\ 0.0]^T\; \LU$ and $\vI(\tin) = [-0.5\ -2.5\ 0.0]^T \; \LU/\TU$. To land successfully, the vehicle must shed significant horizontal momentum while ensuring an upright final attitude. 

Figure~\ref{fig:exI_trjs} shows the converged trajectory for each of the three cases. To reduce clutter, only~$11$ of the~$20$ temporal nodes are shown. The insets~(i) and~(ii) represent free body diagrams of the forces acting on the vehicle at the same two temporal nodes for each case. Figure~\ref{fig:exI_nostc_nodrag} shows the trajectory for the case without aerodynamic effects. From the insets, it is clear that the maneuver resembles a reverse gravity turn, with the thrust pointed nearly anti-parallel to the velocity vector. 

Figure~\ref{fig:exI_nostc_sphere} shows the trajectory for the case with the spherical aerodynamic model. The corresponding insets show aerodynamic forces that are anti-parallel to the velocity vector, and may therefore be interpreted as drag forces. The thrust directions, however, remain consistent with the no-aerodynamic case. Thus, the trajectory in Figure~\ref{fig:exI_nostc_sphere} may be interpreted as the atmospheric counterpart to the reverse gravity turn observed in case (B).

Figure~\ref{fig:exI_nostc_elip} shows the trajectory for the case with the ellipsoidal aerodynamic model. In this case, the aerodynamic force vectors are seen to have a component orthogonal to the velocity vector, and may therefore be interpreted as a composition of a drag and lift force (recall Figure~\ref{fig:aerodynamic_model}). The vehicle is observed to exploit the lift force to bend the trajectory downwards by adjusting its attitude to control the angle of attack. As a consequence, the thrust is no longer aligned with the velocity as it was in the previous two cases, and is instead gimbaled such that the vehicle remains trimmed at an angle of attack that applies the lift force in a desirable direction.

%%%%%%%%%%%%%%%%%%%%%%%%%%%%%%%%%%%%%%%%%%%%%%%%%%%%%%%%%%%%%%%%%%%%%%%%%%%%%%%%%%

\subsection{State-Triggered Constraints Case Study} \label{sec4:exII}

\begin{figure}[t!]
  \begin{subfigure}{0.5\textwidth}
    \input{tikz_stc_III.tex}
    \caption{Case~(B)+(EA): Trajectory}
    \label{fig:stc_noaoa}
  \end{subfigure} \hfil
  \begin{subfigure}{0.5\textwidth}
    \input{tikz_stc_IV.tex}
    \caption{Case~(B)+(EA)+(ST): Trajectory}
    \label{fig:stc_aoa}
  \end{subfigure}
  \begin{subfigure}{0.5\textwidth}
    \centering
     \includegraphics[width=\textwidth]{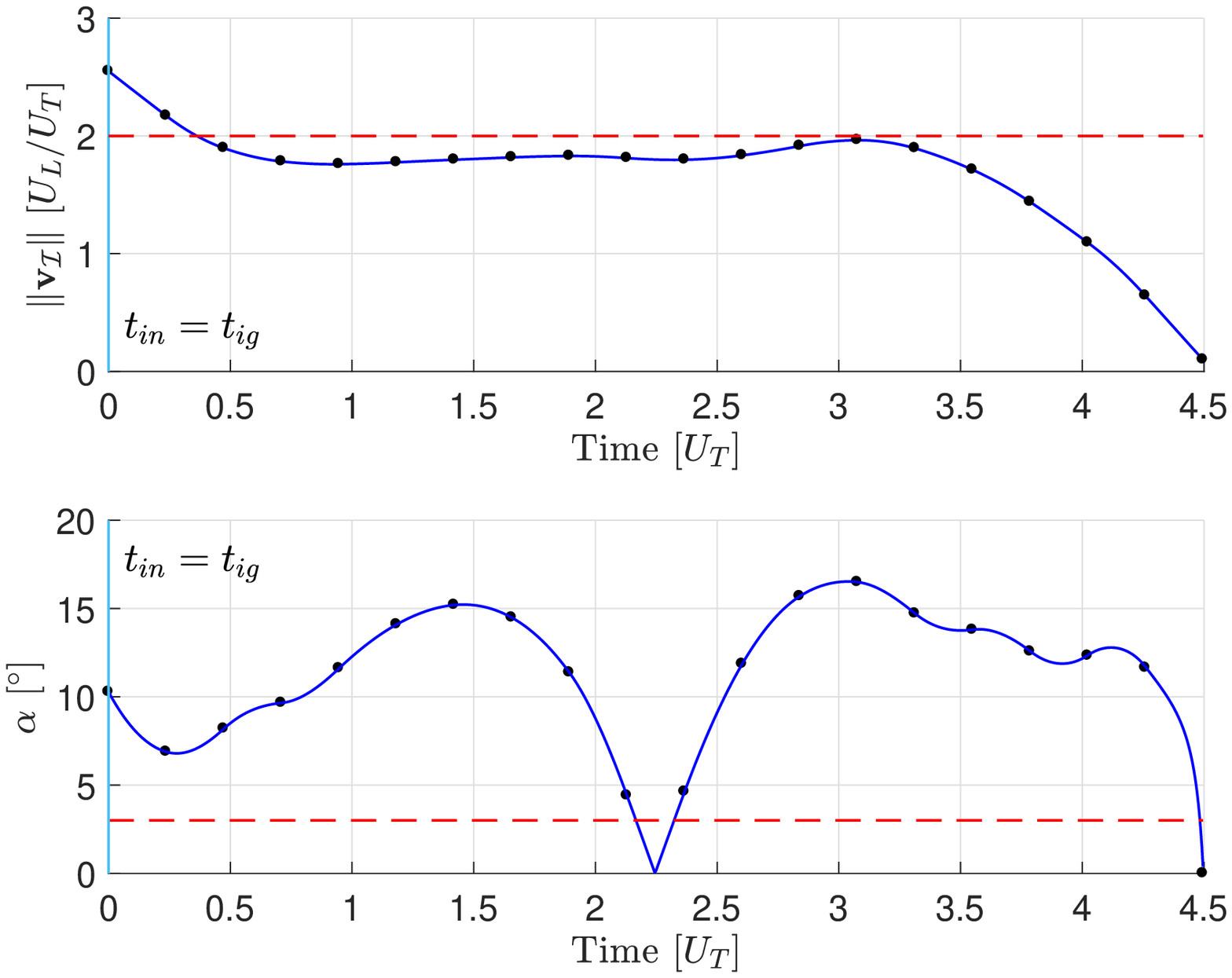}
     \caption{Case~(B)+(EA): STC Parameters}
     \label{fig:exII_nostc_nofit_vandaoa}
  \end{subfigure} \hfil
  \begin{subfigure}{0.5\textwidth}
    \centering
     \includegraphics[width=\textwidth]{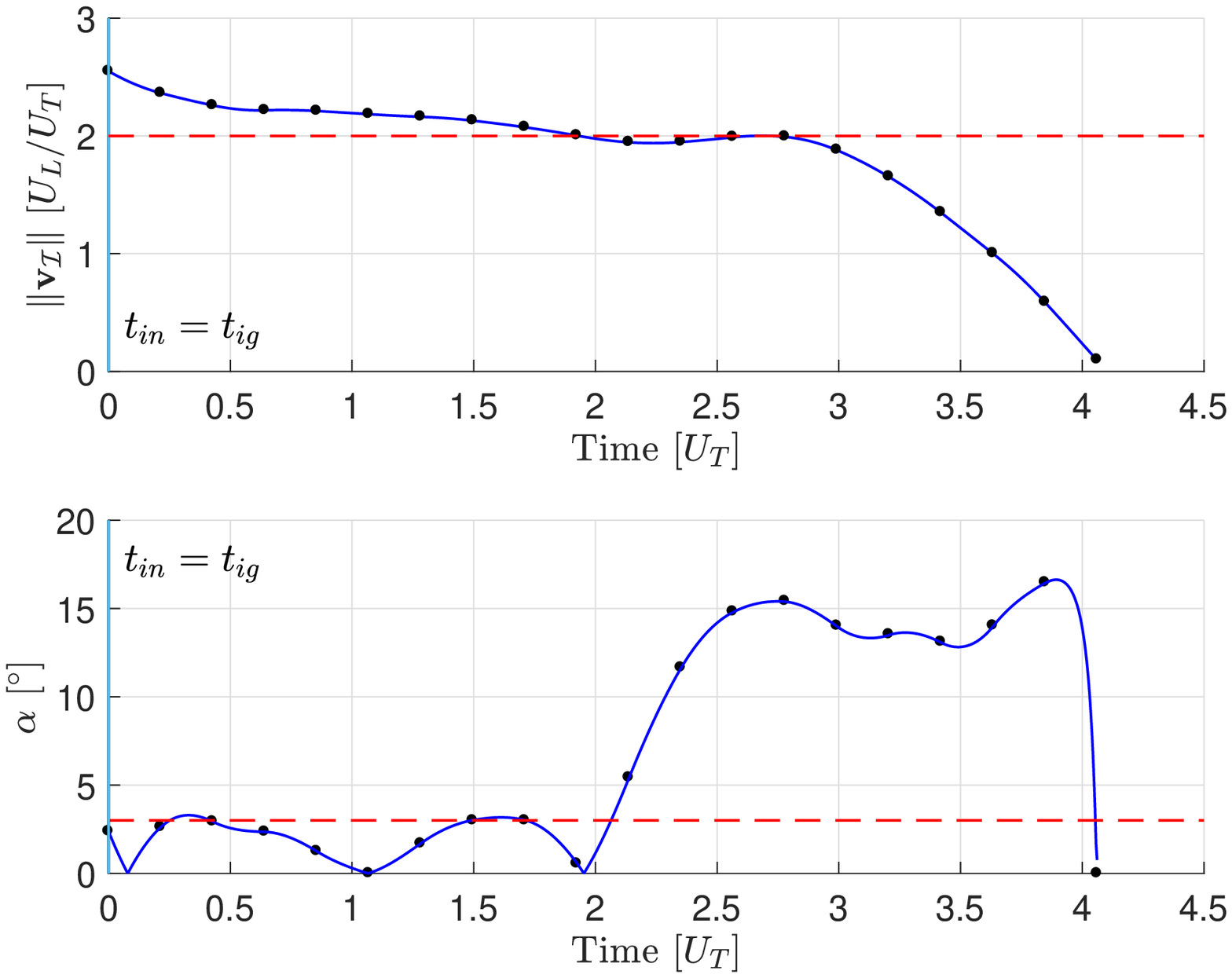}
     \caption{Case~(B)+(EA)+(ST): STC Parameters}
     \label{fig:exII_stc_nofit_vandaoa}
  \end{subfigure}
  \caption{State-Triggered Constraint Case Study: Trajectories and speed-angle of attack histories for cases (B)+(EA) and (B)+(EA)+(ST). This case study uses $\rI(\tin)=[5.33\ 4.5\ 0.0]^T$~$\LU$ and $\vI(\tin)=[-0.5\ -2.5\ 0.0]^T$~$\LU/\TU$. Case (B)+(EA)+(ST) implements the \qalpha~state-triggered constraint introduced in~\eqref{eq:stc_ex_aoa_orig}, which limits the angle of attack to $3.0\dg$ when the speed is greater than $2.0$~$\LU/\TU$. The horizontal dashed red lines in (c) and (d) represent the velocity limit~$\Vaoa$ in the top set of axes, and the angle of attack limit~$\aoamax$ in the bottom set of axes. Note that that these limits are not enforced in (c), and are shown only for reference. Refer to the legend in Figure~\ref{fig:exI_trjs} for additional definitions.}
  \label{fig:exII_trjs}
\end{figure}

This case study highlights the effects of the \qalpha~state-triggered constraint introduced in~\eqref{eq:stc_ex_aoa_orig} on a landing scenario with non-negligible atmospheric effects. We consider the cases (B)+(EA) and (B)+(EA)+(ST), with initial position and velocity vectors $\rI(\tin)=[5.33\ 4.5\ 0.0]^T$~$\LU$ and $\vI(\tin)=[-0.5\ -2.5\ 0.0]^T$~$\LU/\TU$, respectively.

The resulting trajectories are shown in Figures~\ref{fig:stc_noaoa}-\ref{fig:stc_aoa}, with the corresponding speed and angle of attack profiles provided in Figures~\ref{fig:exII_nostc_nofit_vandaoa}-\ref{fig:exII_stc_nofit_vandaoa}. In Figures~\ref{fig:stc_noaoa}-\ref{fig:stc_aoa}, inset~(i) shows the initial condition where the speed is greater than~$\Vaoa$, while inset~(ii) shows the last temporal node at which the STC is active in case (B)+(EA)+(ST). The lift-to-drag ratios displayed in inset~(i) of Figures~\ref{fig:stc_noaoa}-\ref{fig:stc_aoa} indicate that the aerodynamic loading in case (B)+(EA)+(ST) is reduced when the compared to case~(B)+(EA) due to the inclusion of the \qalpha~constraint. The same effect is apparent to a lesser extent in inset~(ii), where the STCs trigger condition is still satisfied for case (B)+(EA)+(ST). In Figure~\ref{fig:stc_aoa}, the temporal nodes after inset~(ii) clearly exhibit larger angles of attack, indicating that the vehicle's speed has dropped below the trigger limit and that the \qalpha~constraint has been disabled.

Figures~\ref{fig:exII_nostc_nofit_vandaoa}-\ref{fig:exII_stc_nofit_vandaoa} show the corresponding speed and angle of attack time histories. Since case (B)+(EA) does not implement the \qalpha~STC, the angle of attack is seen to violate the~$3.0\dg$ limit when the speed is greater than~$2.0$~$\LU/\TU$. In contrast, the angle of attack in case (B)+(EA)+(ST) remains below the prescribed~$3.0\dg$ limit until~$t\approx 2.0$~$\TU$, where the speed drops below the prescribed~$2.0$~$\LU/\TU$ speed limit.

Lastly, notice that Figure~\ref{fig:exII_stc_nofit_vandaoa} shows that case (B)+(EA)+(ST) rides the speed constraint for $t\in\intee{2.0}{3.0}$. Although the trajectory may gain optimality by traveling faster over this time interval, the angle of attack is above the specified $3.0\dg$ limit over this time interval. Consequently, we observe the enforcement of the \textit{contrapositive} of~\eqref{eq:stc_ex_aoa_orig} (i.e. an angle of attack greater than $3.0\dg$ implies that the speed must be less than $2.0$~$\LU/\TU$). Further, certain scenarios may allow the angle of attack to drop below the prescribed limit towards the end of the trajectory, thereby allowing the speed to increase above the STC trigger limit. Such an eventuality would cause the optimization algorithm to add a new constrained phase \textit{without a priori input}, per the discussion in Alternative 3 in~\sref{sec2:stc_ex_app}.

%%%%%%%%%%%%%%%%%%%%%%%%%%%%%%%%%%%%%%%%%%%%%%%%%%%%%%%%%%%%%%%%%%%%%%%%%%%%%%%%%%

\subsection{Free-Ignition-Time Case Study} \label{sec4:exIII}

The third case study considers how the free-ignition time modification affects the trajectory by comparing the cases (B)+(EA) and (B)+(EA)+(FI). The initial position and velocity are identical to those used in~\sref{sec4:exII}. Figure~\ref{fig:exIII_trjs} depicts the converged trajectories of both cases. Inset (i) on each figure corresponds to the ignition time epoch $\tig$, while inset (ii) corresponds to a temporal node in the middle of the maneuver. Case (B)+(EA)+(FI) selects a coast time of $\tc = 0.96$~$\TU$ and as a result observes a reduction in both burn time and fuel cost (see~\sref{sec4:timing}).

\begin{figure}[t!]
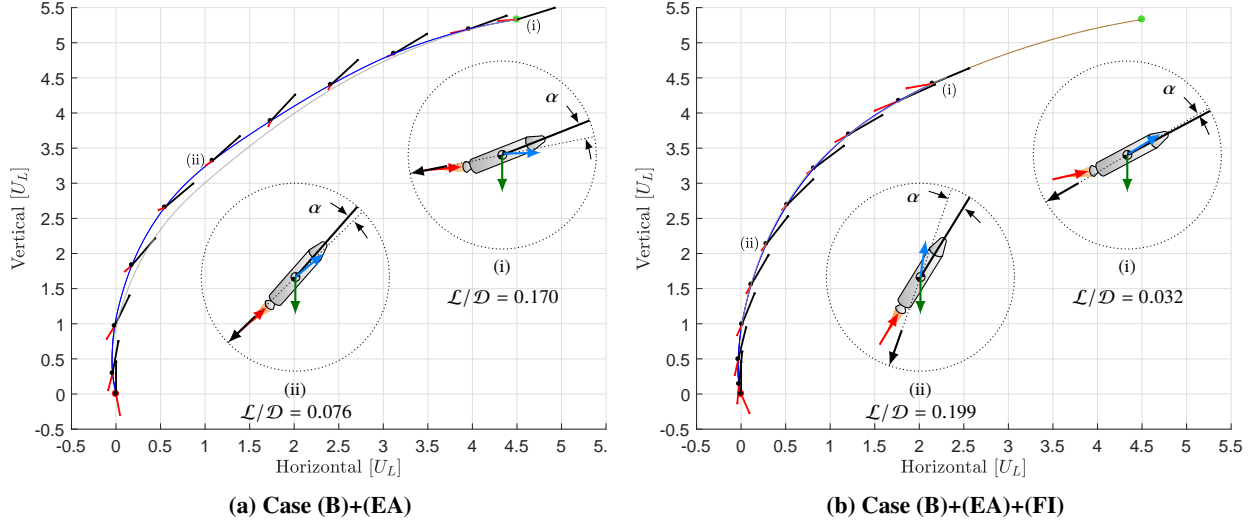

  \begin{subfigure}{0.5\textwidth}
    \input{tikz_stc_III.tex}
    \caption{Case~(B)+(EA)}
    \label{fig:exIII_nostc_nofit}
  \end{subfigure} \hfil
  \begin{subfigure}{0.5\textwidth}
    \input{tikz_stc_V.tex}
    \caption{Case~(B)+(EA)+(FI)}
    \label{fig:exIII_nostc_fit}
  \end{subfigure}
  \caption{Free-Ignition-Time Case Study: Trajectories for case (B)+(EA) and (B)+(EA)+(FI), with $\rI(\tin)=[5.33\ 4.5\ 0.0]^T$~$\LU$ and $\vI(\tin)=[-0.5\ -2.5\ 0.0]^T$~$\LU/\TU$. Case (B)+(EA)+(FI) implements the free-ignition-time modification introduced in \sref{sec2:fit}. Refer to the legend in Figure~\ref{fig:exI_trjs} for additional definitions.}
  \label{fig:exIII_trjs}
\end{figure}

%%%%%%%%%%%%%%%%%%%%%%%%%%%%%%%%%%%%%%%%%%%%%%%%%%%%%%%%%%%%%%%%%%%%%%%%%%%%%%%%%%

\subsection{Trajectory and Computational Performance} \label{sec4:timing}

The trajectory and computational performance data presented in this section were generated for the cases listed in the leftmost column of Table~\ref{tab:performance}. The results were obtained by executing a batch of 10 runs for each case using the initial position and velocity vectors~$\rI(\tin)=[5.33\ 4.5\ 0.0]^T$~$\LU$ and~$\vI(\tin)=[-0.5\ -2.5\ 0.25]^T$~$\LU/\TU$. These initial conditions generate three-dimensional out-of-plane trajectories that render the problem less sparse and thus more computationally intensive than their their planar counterparts presented in~\sref{sec4:exI}-\sref{sec4:exIII}.

Trajectory performance metrics are given in Table~\ref{tab:performance}. For each case, the entries represent the median values of the metrics generated in the batch. Due to the determinism of the proposed algorithm, the standard deviations of the trajectory performance metrics was zero.

To validate each solution against the nonlinear dynamics, we integrate the nonlinear equations of motion using a piecewise linear interpolation of the controls, given in~\eqref{eq:foh_def}. The errors between the integrated trajectory and the discrete states generated by the optimization process are then used as a measure of the solution's feasibility. The position and attitude errors are defined as~$e_{\textit{pos}} \definedas \max_{k \in K}\;\| \rI(t_k) - \rIk \|_2$ and~$e_{\textit{att}} \definedas \max_{k \in K} \; 2 \cos^{-1} \big[\qIB^*(t_k) \otimes \bvec{q}_{\body \leftarrow \inertial,k} \big]$, where~$\bvec{q}_{\body \leftarrow \inertial,k}$ and~$\rIk$ are the discrete solution values, while~$\qIB(t_k)$ and~$\rI(t_k)$ are the corresponding integrated values. 

The computational performance results for the solve step is given in Table~\ref{tab:timing}, and were generated on a 2014 MacBook Pro with a 2.2~GHz Intel Core~i7 processor and~16~GB of RAM. The propagation step was implemented in C++ using the Eigen matrix library~\cite{eigenweb}, and was omitted from Table~\ref{tab:timing} since the maximum propagation time per run was on the order of~$10$ milliseconds. The solve time results were generated in MATLAB using ECOS~\cite{Domahidi2013ecos} and CVX~\cite{cvx}. For each case, the solve times were obtained by totaling the fifth argument reported by the \texttt{cvx\_toc} function over all ECOS calls in a run, and computing the statistics over the entire batch.

\begin{table}[tbh]
  \centering
  \caption{Trajectory Performance Results for Combinations of Problem Features}
  \begin{tabular}{lccccc}
    \hhline{======}
    \multirow{2}*{Case} & Burn         & Fuel          & Successive & $e_{\textit{pos}}$ & $e_{\textit{att}}$ \\
                     & Time ($\TU$) & Consumed (\%) & Iterations & ($\LU$)            & ($\dg$)              \\
    \hline
    (B)                             & $4.42$ & $3.17$ & $4$ & $1.1\text{e-}3$ & $8.0\text{e-}2$     \\
    (B) + (SA)                      & $4.27$ & $3.72$ & $4$ & $2.2\text{e-}2$ & $1.3\text{e-}1$     \\
    (B) + (EA)                      & $4.67$ & $3.55$ & $4$ & $7.8\text{e-}2$ & $1.0\text{e-}0$     \\ % (FI) A1
    (B) + (EA) + (ST)               & $4.06$ & $3.40$ & $7$ & $3.3\text{e-}3$ & $3.7\text{e-}1$     \\ % (FI) B1
    (B) + (ST) + (FI)               & $3.10$ & $3.03$ & $7$ & $3.2\text{e-}5$ & $1.1\text{e-}2$     \\
    (B) + (EA) + (FI)               & $2.56$ & $3.13$ & $5$ & $1.4\text{e-}2$ & $\leq1.0\text{e-}6$ \\ % (FI) A2 ... A1 > A2 (GOOD: (FI) improved optimality)
    (B) + (EA) + (3I)               & $4.89$ & $3.61$ & $4$ & $2.7\text{e-}2$ & $3.2\text{e-}1$     \\ % (FI) C1
    (B) + (EA) + (3I) + (FI)        & $4.41$ & $3.58$ & $4$ & $2.2\text{e-}2$ & $1.1\text{e-}1$     \\ % (FI) C2 ... C1 > C2 (GOOD: (FI) improved optimality)
    (B) + (EA) + (ST) + (FI)        & $2.69$ & $3.21$ & $7$ & $2.1\text{e-}3$ & $\leq1.0\text{e-}6$ \\ % (FI) B2 ... B1 > B2 (GOOD: (FI) improved optimality)
    (B) + (EA) + (ST) + (FI) + (3I) & $3.21$ & $3.27$ & $7$ & $8.0\text{e-}3$ & $\leq1.0\text{e-}6$ \\
    \hhline{======}
  \end{tabular}
  \label{tab:performance}
\end{table}

\begin{table}[tbh]
  \centering
  \caption{Computational Performance Results Averaged Over 10 Runs}
  \begin{tabular}{lccccc}
  	\hhline{======}
    \multirow{2}*{Case}             & \multicolumn{4}{c}{Solve Time [s]} \\ 
                                    & Min & Max & Median & Std \\
    \hline
    (B)                             & $0.164$ & $0.287$ & $0.164$ & $0.007$ \\
    (B) + (SA)                      & $0.170$ & $0.210$ & $0.174$ & $0.012$ \\
    (B) + (EA)                      & $0.271$ & $0.312$ & $0.276$ & $0.012$ \\ % (3I) A1
    (B) + (EA) + (ST)               & $0.496$ & $0.565$ & $0.504$ & $0.020$ \\
    (B) + (ST) + (FI)               & $0.643$ & $0.682$ & $0.654$ & $0.011$ \\
    (B) + (EA) + (FI)               & $0.233$ & $0.263$ & $0.241$ & $0.008$ \\ % (3I) B1
    (B) + (EA) + (3I)               & $0.299$ & $0.329$ & $0.312$ & $0.011$ \\ % (3I) A2 ... A1 < A2 (REALLY BAD: (3I) increased time, reduced optimality)
    (B) + (EA) + (3I) + (FI)        & $0.224$ & $0.253$ & $0.229$ & $0.008$ \\ % (3I) B2 ... B1 > B2 (BAD:        (3I) decreased time, reduced optimality)
    (B) + (EA) + (ST) + (FI)        & $0.494$ & $0.537$ & $0.505$ & $0.012$ \\ % (3I) C1
    (B) + (EA) + (ST) + (FI) + (3I) & $0.602$ & $0.632$ & $0.618$ & $0.009$ \\ % (3I) C2 ... C1 < C2 (REALLY BAD: (3I) increased time, reduced optimality)
  	\hhline{======}
  \end{tabular}
  \label{tab:timing}
\end{table}

We conclude with three final observations. First, in the cases presented, the inclusion of feature (ST) yielded a decrease in both burn time and fuel consumption. This result is counter intuitive since one would expect the optimal cost of a problem to remain the same or increase when additional constraints are added. Since Problem~\ref{prob:ncvx} may have multiple local minima, we posit that the inclusion of the \qalpha~state-triggered constraint forced the iterative solution process towards a more optimal local minima. On the other hand, the inclusion of feature (FI) reduced the fuel cost. This result agrees with intuition since the free-ignition-time modification effectively adds a degree of freedom to the problem.

Second, we note that the 3-DoF initialization approach (i.e. the inclusion of feature (3I)) did not yield a clear improvement in optimality or computational performance. In fact, cases (B)+(EA)+(3I), (B)+(EA)+(3I)+(FI), and (B)+(EA)+(ST)+(FI)+(3I) all resulted in less optimal trajectories, while cases (B)+(EA)+(3I) and (B)+(EA)+(ST) +(FI)+(3I) also increased the solve time. We conclude that the straight-line initialization approach initializes the algorithm in a more favorable region of attraction, but stress that this may not be the case for different scenarios.

Lastly, the timing results presented in Table~\ref{tab:timing} show a maximum solve time of~$0.7$ seconds and a standard deviation on the order of milliseconds and were all obtained using the same algorithm parameters (e.g.~$\wvc$, $\Wtr$, $\KK$). We argue that these results are an important step in demonstrating the efficacy of the successive convexification methodology for real time autonomous applications. Ultimately, results obtained on representative flight hardware will be crucial in accurately assessing the viability of the proposed methodology for on-board computation in real-world applications.

%%%%%%%%%%%%%%%%%%%%%%%%%%%%%%%%%%%%%%%%%%%%%%%%%%%%%%%%%%%%%%%%%%%%%%%%%%%%%%%%%%%%%%%
%%%%%%%%%%%%%%%%%%%%%%%%%%%%%%%%%%%%%%%%%%%%%%%%%%%%%%%%%%%%%%%%%%%%%%%%%%%%%%%%%%%%%%%
%%%%%%%%%%%%%%%%%%%%%%%%%%%%%%%%%%%%%%%%%%%%%%%%%%%%%%%%%%%%%%%%%%%%%%%%%%%%%%%%%%%%%%%

\section{Concluding Remarks} \label{sec:conclusion}

This paper presents a real-time successive convexification algorithm for a generalized free-final-time 6-DoF powered descent guidance problem, and introduces three primary contributions: (i) a free-ignition-time modification that allows the algorithm to determine when to begin the burn phase, (ii) an ellipsoidal aerodynamics model that provides a computationally tractable way to model lift and drag forces, and (iii) a continuous formulation for state-triggered constraints. Contribution (iii) allows continuous optimization problems to be formulated using conditionally enforced constraints, and was motivated by landing scenarios that necessitate velocity-triggered angle of attack and range-triggered line of sight constraints.

Three simulation case studies are presented, each illustrating one of the primary contributions of this paper. The corresponding trajectory and computational performance results show that the proposed algorithm can successfully compute trajectories in under~$0.7$~seconds for the problem features considered. While additional work is required to provide convergence guarantees and to quantify the optimality of the computed trajectories, we argue that our results demonstrate the efficacy of the successive convexification approach for real-time powered descent guidance applications.

%%%%%%%%%%%%%%%%%%%%%%%%%%%%%%%%%%%%%%%%%%%%%%%%%%%%%%%%%%%%%%%%%%%%%%%%%%%%%%%%%%%%%%%
%%%%%%%%%%%%%%%%%%%%%%%%%%%%%%%%%%%%%%%%%%%%%%%%%%%%%%%%%%%%%%%%%%%%%%%%%%%%%%%%%%%%%%%
%%%%%%%%%%%%%%%%%%%%%%%%%%%%%%%%%%%%%%%%%%%%%%%%%%%%%%%%%%%%%%%%%%%%%%%%%%%%%%%%%%%%%%%

\bibliography{bibliography}

%%%%%%%%%%%%%%%%%%%%%%%%%%%%%%%%%%%%%%%%%%%%%%%%%%%%%%%%%%%%%%%%%%%%%%%%%%%%%%%%%%%%%%%
%%%%%%%%%%%%%%%%%%%%%%%%%%%%%%%%%%%%%%%%%%%%%%%%%%%%%%%%%%%%%%%%%%%%%%%%%%%%%%%%%%%%%%%
%%%%%%%%%%%%%%%%%%%%%%%%%%%%%%%%%%%%%%%%%%%%%%%%%%%%%%%%%%%%%%%%%%%%%%%%%%%%%%%%%%%%%%%

\end{document}